\documentclass[nohyperref]{article}

\usepackage{microtype}
\usepackage{graphicx}
\usepackage{booktabs} 

\usepackage{hyperref}

\usepackage[accepted]{icml2022}

\usepackage{amsmath}
\usepackage{amssymb}
\usepackage{mathtools}
\usepackage{amsthm}
\usepackage{multirow}
\usepackage{subcaption}

\usepackage[capitalize,noabbrev]{cleveref}

\theoremstyle{plain}
\newtheorem{theorem}{Theorem}[section]
\newtheorem{proposition}[theorem]{Proposition}
\newtheorem{lemma}[theorem]{Lemma}

\theoremstyle{definition}
\newtheorem{definition}[theorem]{Definition}
\newtheorem{assumption}[theorem]{Assumption}
\theoremstyle{remark}
\newtheorem{remark}[theorem]{Remark}

\usepackage[textsize=tiny]{todonotes}

\icmltitlerunning{Gradient-Free Method for Heavily Constrained Nonconvex Optimization}

\begin{document}

\twocolumn[
\icmltitle{Gradient-Free Method for Heavily Constrained Nonconvex Optimization}



\icmlsetsymbol{equal}{*}

\begin{icmlauthorlist}
\icmlauthor{Wanli Shi}{yyy,comp}
\icmlauthor{Hongchang Gao}{sch}
\icmlauthor{Bin Gu}{yyy,comp}
\end{icmlauthorlist}

\icmlaffiliation{yyy}{Nanjing University of Information Science and Technology, Jiangsu, China}
\icmlaffiliation{comp}{MBZUAI, Abu Dhabi, UAE}
\icmlaffiliation{sch}{Department of Computer and Information Sciences, Temple University, PA, USA}

\icmlcorrespondingauthor{Bin Gu}{jsgubin@gmail.com}

\icmlkeywords{Constrained Optimization, Zeroth-order Gradient, Nonconvex Optimization}

\vskip 0.3in
]



\printAffiliationsAndNotice{}  

\begin{abstract}
Zeroth-order (ZO) method has been shown to be a powerful method for solving the optimization problem where explicit expression of the gradients is difficult or infeasible to obtain. Recently, due to the practical value of the constrained problems, a lot of ZO Frank-Wolfe or projected ZO methods have been proposed. However, in many applications, we may have a very large number of nonconvex white/black-box constraints, which makes the existing zeroth-order methods extremely inefficient (or even not working) since they need to inquire function value of all the constraints and project the solution to the complicated feasible set. In this paper, to solve the nonconvex problem with a large number of white/black-box constraints, we proposed a doubly stochastic zeroth-order gradient method (DSZOG) with momentum method and adaptive step size. Theoretically, we prove DSZOG can converge to the $\epsilon$-stationary point of the constrained problem. Experimental results in two applications demonstrate the superiority of our method in terms of training time and accuracy compared with other ZO methods for the constrained problem.
\end{abstract}

\section{Introduction}

Zeroth-order (gradient-free) method is a powerful method for solving the optimization problem where explicit expression of the gradients are difficult or infeasible to obtain, such as bandit feedback analysis \cite{agarwal2010optimal}, reinforcement learning \cite{choromanski2018structured}, and adversarial attacks on black-box deep neural networks \cite{chen2017zoo,liu2018zeroth_var}. Recently, more and more zeroth-order gradient algorithms have been proposed and achieved great success, such as  \cite{ghadimi2013stochastic,wang2018stochastic,gu2016zeroth,liu2018zeroth_var,huang2020accelerated,gu2021black,wei2021accelerated,gu2021optimizing}. 

Due to several motivating applications, the study of the zeroth-order methods in constrained optimization has gained great attention.  For example, ZOSCGD\citep{balasubramanian2018zeroth} uses the zeroth-order method to approximate the unbiased stochastic gradient of the objective, and then uses the Frank-Wolfe framework to update the parameters. \cite{gao2020can,huang2020accelerated2} apply the variance reduction technique \cite{fang2018spider,nguyen2017sarah} or momentum method in ZOSCGD and obtain a better convergence performance. In addition, ZOSPGD \cite{liu2018zeroth_pro} uses the zeroth-order gradient to update the parameters and then projects the parameters onto the feasible subset. The variance reduction and momentum methods are also used to obtain a better performance \cite{huang2020accelerated}.
We have summarized several representative zeroth-order methods for constrained optimization in Table \ref{tab:zo_cop}.

\begin{table*}[]
	\small
	\centering
	\setlength{\tabcolsep}{4.2mm}
	\caption{Representative zeroth order methods for constrained optimization problems, where N/C means nonconvex/convex, W/B means white/black-box function, and the last column shows the size of the constraints.}
	\label{tab:zo_cop}
	\begin{tabular}{c|l|c|c|c|c}
		\hline
		Framework                    & Algorthm                                                               & Reference                                                      & Objective                  & Constraints                                                                  &  Size    \\ \hline
		\multirow{6}{*}{Frank-Wolfe} & ZOSCGD                                                                  & \cite{balasubramanian2018zeroth}              & N/C                  & \begin{tabular}[c]{@{}c@{}}C\\ W\end{tabular}                      & Small                  \\ \cline{2-6} 
		& FZFW                                                                    & \multirow{3}{*}{\cite{gao2020can}}            & \multirow{3}{*}{N/C} & \multirow{3}{*}{\begin{tabular}[c]{@{}c@{}}C\\ W\end{tabular}}     & \multirow{3}{*}{Small} \\ \cline{2-2}
		& FZCGS                                                                   &                                                                &                          &                                                                              &                        \\ \cline{2-2}
		& FCGS                                                                    &                                                                &                          &                                                                              &                        \\ \cline{2-6} 
		& Acc-SZOFW                                                               & \multirow{2}{*}{\cite{huang2020accelerated2}} & \multirow{2}{*}{N/C} & \multirow{2}{*}{\begin{tabular}[c]{@{}c@{}}C\\ W\end{tabular}}     & \multirow{2}{*}{Small} \\ \cline{2-2}
		& Acc-SZOFW*                                                              &                                                                &                          &                                                                              &                        \\ \hline
		\multirow{2}{*}{Projected}   & ZOPSGD                                                                   & \cite{liu2018zeroth_pro}                     & N/C                  & \begin{tabular}[c]{@{}c@{}}C\\ W\end{tabular}                      & Small                  \\ \cline{2-6} 
		& AccZOMDA                                                          & \cite{huang2020accelerated}                            & N/C                  & \begin{tabular}[c]{@{}c@{}}C\\ W\end{tabular}                      & Small                  \\ \hline
		\multirow{2}{*}{Penalty}     & \multirow{2}{*}{DSZOG} & \multirow{2}{*}{Ours}                                              & \multirow{2}{*}{N/C} & \multirow{2}{*}{\begin{tabular}[c]{@{}c@{}}N/C\\ W/B\end{tabular}} & \multirow{2}{*}{Large} \\
		&                                                                         &                                                                &                          &                                                                              &                        \\ \hline
	\end{tabular}
\end{table*}
However, all these methods are not scalable for the problems with a large number of constraints. On the one hand, they all need to evaluate the values of all the constraints in each iteration. On the other hand, the projected gradient methods and the Frank-Wolfe methods need to solve a subproblem in each iteration. These makes the existing methods time-consuming to find a feasible point. What's worse, all these methods need the constraints to be convex white-box functions. However, in many real-world applications, the constraints could be nonconvex or black-box functions, which means existing methods are extremely inefficient or even not working. Therefore, how to effectively solve the nonconvex constrained problem with a large number of nonconvex/convex white/black-box constraints, which is denoted as heavily constrained problem, by using the ZO method is still an open problem.

In this paper, to solve the heavily constrained nonconvex optimization problem efficiently, we propose a new ZO algorithms called doubly stochastic zeroth-order gradient method (DSZOG). Specifically, we give a probability distribution over all the constraints and rewrite the original problem as a nonconvex-strongly-concave minimax problem \cite{lin2020gradient,wang2020zeroth,huang2020accelerated,guo2021stochastic} by using the penalty method. We sample a batch of training points uniformly and a batch of constraints according to the distribution to calculate the zeroth-order gradient of the penalty function w.r.t model parameters and then sample a batch of constraints uniformly to calculate the stochastic gradient of penalty function w.r.t the probability distribution. Then, gradient descent and projected gradient ascent can be used to update model parameters and probability distribution, respectively. In addition, we also use the exponential moving average (EMA) method and adaptive stepsize \cite{guo2021stochastic,huang2020accelerated}, which benefits our method from the variance reduction and adaptive convergence. Theoretically, we prove DSZOG can converge to the $\epsilon$-stationary point of the constrained problem. Experimental results in two applications demonstrate the superiority of our method in terms of training time and accuracy compared with other ZO methods for constrained problem.

\noindent \textbf{Contributions.} 
We summarized the main contributions of this paper as follows:
\begin{enumerate}
	\item We propose a doubly stochastic zeroth-order gradient method to solve the heavily constrained nonconvex problem. By introducing a stochastic layer into the constraints, our method is scalable and efficient for the heavily constrained nonconvex problem. 
	\item By using the exponential moving average method and adaptive stepsize, our method enjoys the benefits of variance reduction and adaptive convergence.
	\item We prove DSZOG can converge to the $\epsilon$-stationary point of the constrained problem. Experimental results also demonstrate the superiority of our methods in terms of accuracy and training time.
\end{enumerate}

\section{Related Works}
\subsection{Zeroth-Order Methods}
Zeroth-order methods are powerful methods to solve several machine learning problems, where the explicit gradients are difficult or infeasible to obtain. Based on Gaussian smoothing method, \cite{ghadimi2013stochastic,duchi2015optimal,nesterov2017random,gu2016zeroth,gu2021black,wei2021accelerated,gu2021optimizing} propose several zeroth-order method which only needs function values to estimate the gradients. To deal with nonsmooth optimization problem, some zeroth-order proximal gradient methods \cite{ghadimi2016mini,ji2019improved} and ADMM methods \cite{gao2018information,liu2018zeroth} have been proposed. In addition, to solve the constrained optimization problems, the zeroth-order projection method \cite{liu2018zeroth_pro} and the zeroth-order Frank-Wolfe methods \cite{balasubramanian2018zeroth,chen2020frank} have been proposed. More recently, based on the variance reduced techniques, some accelerated zeroth-order stochastic methods have been proposed. \cite{gu2021black} proposed a new framework to reduce the query complexities of zeroth-order gradient methods for convex and nonconvex objectives. In addition, this new framework can be used in various ZO method and has a better convergence performance.

\subsection{Variance Reduction and Momentum Methods}
To accelerate stochastic gradient descent method, variance reduction methods such as SAG \cite{roux2012stochastic}, SAGA \cite{defazio2014saga}, SVRG \cite{johnson2013accelerating}, SARAH \cite{nguyen2017sarah} have been proposed. Recently, due to the widely existence of nonconvex optimization, several variance reduction methods for nonconvex optimization have been proposed \cite{allen2017katyusha,lei2017non,fang2018spider,wang2019spiderboost,zhou2018stochastic}. Another method to accelerate the stochastic gradient method is to use momentum-based method. For convex and nonconvex optimization problem, various momentum-based methods have been proposed, e.g. APCG \cite{lin2014accelerated}, Katyusha \cite{allen2017katyusha}, STORM \cite{cutkosky2019momentum}, NIGHT \cite{cutkosky2020momentum}, Hybrid-SGD \cite{tran2021hybrid}, etc. Due to the superiority of variance reduction and momentum methods, they have been widely used in zeroth-order gradient methods and have achieved great success.

\section{Preliminaries}
\subsection{Problem Setting}
In this paper, we consider the following nonconvex constrained problem,
\begin{align}\label{constrained_problem}
&\min_{\boldsymbol{w}} \ f_0(\boldsymbol{w}):=\dfrac{1}{n}\sum_{i=1}^{n}\ell_i(\boldsymbol{w}), \\
& s.t.  \ f_j(\boldsymbol{w})\leq 0, \ j=1,\cdots,m, \nonumber
\end{align}
where $\boldsymbol{w}\in\mathbb{R}^d$ is the optimization variable, $\{\ell_i(\boldsymbol{w})\}_{i=1}^n$ are $n$ component functions. In addition, $f_0:\mathbb{R}^d\mapsto\mathbb{R}$ is a nonconvex and black-box function. $f_j:\mathbb{R}^d\mapsto\mathbb{R}$, $(j=1,\cdots,m)$, is nonconvex/convex and  white/black-box function. We can denote such problem as heavily constrained problem.

\subsection{Reformulate the Constrained Problem}
To solve the constrained problem, the penalty method is one of the main approaches and has achieved great success \cite{clarkson2012sublinear,cotter2016light,shi2021improved}. Specifically, the penalty method reformulates the problem by adding a new term onto the objective to penalize the constraints and the solves the new problem to find a KKT point. Based on the penalty method, we reformulate the constrained optimization problem \ref{constrained_problem} as the following minimax problem over a probability distribution \cite{clarkson2012sublinear,cotter2016light}
\begin{align}\label{mini_opt}
\min_{\boldsymbol{w}}\max_{\boldsymbol{p}\in\Delta^m}\mathcal{L}(\boldsymbol{w},\boldsymbol{p})=&f_0(\boldsymbol{w})+\beta\varphi(\boldsymbol{w},\boldsymbol{p})  -\dfrac{\lambda}{2}\|\boldsymbol{p}\|_2^2,
\end{align} 
where $\beta>0$, $\lambda>0$, $\varphi(\boldsymbol{w},\boldsymbol{p})=\sum_{j=1}^{m}p_j\phi_j(\boldsymbol{w})$, $\phi_j(\boldsymbol{w})=(\max\{f_j(\boldsymbol{w}),0\})^2$ is the penalty function on $f_j$, $\Delta^m:=\{\boldsymbol{p}|\sum_{j=1}^{m}p_j=1, 0\leq p_j\leq 1,\forall j\in[m]\}$ is the $m$-dimensional simplex and $\boldsymbol{p}=[p_1,\cdots,p_m]\in\Delta^m$. We add an additional term 
 term $-\dfrac{1}{2}\|\boldsymbol{p}\|_2^2$ in problem \ref{mini_opt} to ensure $\mathcal{L}$ is strongly concave on $\boldsymbol{p}$. 

Since we can only obtain the values of the objective and constraints, stochastic zeroth-order gradient method is one of the effective ways to solve this problem. However, calculating the stochastic zeroth-order gradient of $\mathcal{L}$ needs to inquire the function values of all the constraints, which has a very high time complexity if $m$ is very large. This make it time-consuming.

\section{Proposed Method}
\subsection{Doubly Stochastic Zeroth-order Gradient Method}

To solve problem \ref{mini_opt} efficiently, we introduce the another stochastic layer to the constraints. Specifically, since the minimax problem \ref{mini_opt} contains two finite sums, i.e., $f_0(\boldsymbol{w})=1/n\sum_{i=1}^n\ell_i(\boldsymbol{w})$ and $\varphi(\boldsymbol{w},\boldsymbol{p})=\sum_{j=1}^mp_j\phi_j(\boldsymbol{w})$, we can calculate their stochastic zeroth-order gradient, respectively, and then combine these two gradient to obtain the stochastic zeroth-order gradient of $\mathcal{L}$. 

We can calculate the stochastic zeroth-order gradient of $f_0(\boldsymbol{w})$ and $\varphi(\boldsymbol{w},\boldsymbol{p})$ as follows,
\begin{align}
&G_{\mu}^{f}(\boldsymbol{w}_{t},\ell_i,\boldsymbol{u})=\dfrac{\ell_i(\boldsymbol{w}_{t}+\mu\boldsymbol{u})-\ell_i(\boldsymbol{w}_{t})}{\mu}\boldsymbol{u}, \\ 
&G_{\mu}^{\varphi}(\boldsymbol{w}_{t},\boldsymbol{p},f_j,\boldsymbol{u})=\dfrac{\phi_j(\boldsymbol{w}_{t}+\mu\boldsymbol{u})-\phi_j(\boldsymbol{w}_{t})}{\mu}\boldsymbol{u},
\end{align} 
by sampling $\ell_i$ uniformly, and $f_j$ according to $\boldsymbol{p}$, where $\mu>0$ and $\boldsymbol{u}\sim\mathcal{N}(0,\boldsymbol{1}_d)$. Then, combining these two terms, we can obtain the stochastic zeroth-order gradient of $\mathcal{L}$ w.r.t. $\boldsymbol{w}$ as follows,
\begin{align}
&G_{\mu}^{\mathcal{L}}(\boldsymbol{w}_t,\boldsymbol{p}_t,\ell_i,f_j, \boldsymbol{u})\nonumber\\
=&G_{\mu}^{f}(\boldsymbol{w}_t,\ell_i,\boldsymbol{u})+\beta G_{\mu}^{\varphi}(\boldsymbol{w}_t,\boldsymbol{p}_t,f_j,\boldsymbol{u}).
\end{align}
To reduce the variance, we can sample a batch of $\ell_i$, $f_j$ and $\boldsymbol{u}_k$ to calculate the zeroth-order gradient. Given $q>0$, $\mathcal{M}_1\subseteq[n]$ and $\mathcal{M}_2\sim\boldsymbol{p}\subseteq[m]$, we have 
\begin{align}
&G_{\mu}^{\mathcal{L}}(\boldsymbol{w}_{t},\boldsymbol{p}_{t},\ell_{\mathcal{M}_1},f_{\mathcal{M}_2},\boldsymbol{u}_{[q]})\nonumber\\
=&\dfrac{1}{q|\mathcal{M}_1|}\sum_{i\in\mathcal{M}_1}\sum_{k=1}^{q}G_{\mu}^{f}(\boldsymbol{w}_{t},\ell_i,\boldsymbol{u}_k)\nonumber\\
&+\dfrac{\beta}{q|\mathcal{M}_2|}\sum_{j\in\mathcal{M}_2}\sum_{k=1}^{q}G_{\mu}^{\varphi}(\boldsymbol{w}_{t},\boldsymbol{p}_{t},f_j,\boldsymbol{u}_k),
\end{align} 
Then, the gradient descent can be used to update $\boldsymbol{w}$ by using the following rule 
\begin{align}
\boldsymbol{w}_{t+1}=\boldsymbol{w}_{t}-\eta_{\boldsymbol{w}}	G_{\mu}^{\mathcal{L}}(\boldsymbol{w}_t,\boldsymbol{p}_t,\ell_{\mathcal{M}_1},f_{\mathcal{M}_2},\boldsymbol{u}_{[q]}).
\end{align}

Then, in each iteration, we randomly sample a constraint $f_j(\boldsymbol{w})$ to calculate the stochastic gradient of $\mathcal{L}$ w.r.t. $\boldsymbol{p}$ by using 
\begin{align}
H(\boldsymbol{w}_{t},\boldsymbol{p}_{t},f_j)=\beta m \boldsymbol{e}_j\phi_j(\boldsymbol{w}_t)-\lambda \boldsymbol{p}_{t},
\end{align}
where $\boldsymbol{e}_j$ is the $j$th $m$-dimensional standard unit basis vector. Mini-batch can be also used to reduce variance. Assume we have the randomly sampled index set $\mathcal{M}_3\subseteq[m]$, the mini-batch gradient of $\mathcal{L}$ w.r.t $\boldsymbol{p}$ becomes 
\begin{align}
H(\boldsymbol{w}_{t},\boldsymbol{p}_{t},f_{\mathcal{M}_3})=\dfrac{\beta m}{|\mathcal{M}_3|}\sum_{j\in\mathcal{M}_3} \boldsymbol{e}_j\phi_j(\boldsymbol{w}_t)-\lambda \boldsymbol{p}_{t}.
\end{align}
Then we can perform gradient ascent by using the following rules,
\begin{align}
{\boldsymbol{p}}_{t+1}&=\mathcal{P}_{\Delta^m}(\boldsymbol{p}_t+\eta_{\boldsymbol{p}}H(\boldsymbol{w}_{t},\boldsymbol{p}_t,f_{\mathcal{M}_3})),
\end{align} to update $\boldsymbol{p}$, where $\mathcal{P}_{\Delta^m}(\cdot)$ denotes the projection onto $\Delta^m$ and is easy to calculate.

Note that since $m$ and $n$ are sufficient large in this problem, $G_{\mu}^{\mathcal{L}}(\boldsymbol{w}_t,\boldsymbol{p}_t,\ell_{\mathcal{M}_1},f_{\mathcal{M}_2},\boldsymbol{u}_{[q]})$ and  $H(\boldsymbol{w}_{t},\boldsymbol{p}_{t},f_{\mathcal{M}_3})$ can be viewed as the unbiased estimation of the gradients of $\mathcal{L}$ w.r.t $\boldsymbol{w}$ and $\boldsymbol{p}$, respectively.

\subsection{Momentum and Adaptive Step Size} 
To further improve our method, we use exponential moving average (EMA) method \cite{wang2017stochastic,liu2020adam,cutkosky2020momentum,guo2021stochastic} and adaptive stepsize. 
\begin{algorithm}[!ht]
	\caption{ Doubly Stochastic Zeroth-order Gradient (DSZOG).} 
	\renewcommand{\algorithmicrequire}{\textbf{Input:}}
	\renewcommand{\algorithmicensure}{\textbf{Output:}}
	\begin{algorithmic}[1] 
		\REQUIRE  $T$, $|\mathcal{M}_1|,|\mathcal{M}_2|,|\mathcal{M}_3|$, $\beta \geq1$, $q$, $\mu$, $\lambda=1e-6$, $b\in(0,1)$, $c=1e-8$, $a\in(0,1)$, $\eta_{w}$ and $\eta_{\boldsymbol{p}}$.
		\ENSURE $\boldsymbol{w}_{T}$.
		\STATE Initialize $\boldsymbol{w}_1$.
		\STATE Initialize $\boldsymbol{p}_1=\boldsymbol{p}^*(\boldsymbol{w}_1)$ by solving the strongly concave problem.
		\STATE Initialize $\boldsymbol{z}_{\boldsymbol{w}}^{1}=G_{\mu}^{\mathcal{L}}({\boldsymbol{w}}_{1},\boldsymbol{p}_{1},\ell_{\mathcal{M}_1},f_{\mathcal{M}_2},\boldsymbol{u}_{[q]})$ and $\boldsymbol{z}_{\boldsymbol{p}}^{1}=H(\boldsymbol{w}_{1},{\boldsymbol{p}}_{1},f_{\mathcal{M}_3})$.
		
		\FOR{$t=1,\cdots,T$}
		
		\STATE $\boldsymbol{w}_{t+1}=\boldsymbol{w}_{t}-\eta_{\boldsymbol{w}}\dfrac{\boldsymbol{z}_{\boldsymbol{w}}^t}{\sqrt{\|\boldsymbol{z}_{\boldsymbol{w}}^t\|_2}+c}$.
		
		\STATE $\hat{\boldsymbol{p}}_{t+1}=\mathcal{P}_{\Delta^m}(\boldsymbol{p}_t+\eta_{\boldsymbol{p}}\dfrac{\boldsymbol{z}_{\boldsymbol{p}}^t}{\sqrt{\|\boldsymbol{z}_{\boldsymbol{p}}^t\|_2}+c})$.
		
		
		\STATE $\boldsymbol{p}_{t+1}=(1-a)\boldsymbol{p}_{t}+a\hat{\boldsymbol{p}}_{t+1}$.

		\STATE Randomly sample $\boldsymbol{u}_1,\cdots,\boldsymbol{u}_q\sim\mathcal{N}(0,\boldsymbol{1}_d)$.
		\STATE Randomly sample a index set $\mathcal{M}_1\subseteq[n]$ of $\ell_i$.   
		
		\STATE Sample a constraint index set $\mathcal{M}_2\sim {\boldsymbol{p}}_{t+1}\subseteq[m]$.
		
		\STATE Randomly sample a constraint index set $\mathcal{M}_3$.
		
		\STATE Calculate $G_{\mu}^{\mathcal{L}}(\boldsymbol{w}_{t+1},\boldsymbol{p}_{t+1},\ell_{\mathcal{M}_1},f_{\mathcal{M}_2},\boldsymbol{u}_{[q]})=\dfrac{1}{q|\mathcal{M}_1|}\sum_{i\in\mathcal{M}_1}\sum_{k=1}^{q}G_{\mu}^{f}(\boldsymbol{w}_{t+1},\ell_i,\boldsymbol{u}_k)+\dfrac{\beta}{q|\mathcal{M}_2|}\sum_{j\in\mathcal{M}_2}\sum_{k=1}^{q}G_{\mu}^{\varphi}(\boldsymbol{w}_{t+1},\boldsymbol{p}_{t+1},f_j,\boldsymbol{u}_k).$
		\STATE Calculate $H(\boldsymbol{w}_{t+1},\boldsymbol{p}_{t+1},f_{\mathcal{M}_3})=\dfrac{\beta m}{|\mathcal{M}_3|}\sum_{j\in\mathcal{M}_3} \boldsymbol{e}_j\phi_j(\boldsymbol{w}_{t+1})-\lambda \boldsymbol{p}_{t+1}$.
		
		\STATE $\boldsymbol{z}_{\boldsymbol{w}}^{t+1}=(1-b)\boldsymbol{z}_{\boldsymbol{w}}^t+bG_{\mu}^{\mathcal{L}}({\boldsymbol{w}}_{t+1},\boldsymbol{p}_{t+1},\ell_{\mathcal{M}_1},f_{\mathcal{M}_2},\boldsymbol{u}_{[q]})$.
		\STATE $\boldsymbol{z}_{\boldsymbol{p}}^{t+1}=(1-b)\boldsymbol{z}_{\boldsymbol{p}}^t+bH(\boldsymbol{w}_{t+1},{\boldsymbol{p}}_{t+1},f_{\mathcal{M}_3})$.
		\ENDFOR
	\end{algorithmic}
	\label{alg:DSZOG3}
\end{algorithm}
We use the following exponential moving average (EMA) method on the zeroth-order and first-order gradient to smooth out short-term fluctuations, highlight longer-term trends and reduce the variance of stochastic gradient \cite{wang2017stochastic,guo2021stochastic}
\begin{align}
\boldsymbol{z}_{\boldsymbol{w}}^{t+1}=&(1-b)\boldsymbol{z}_{\boldsymbol{w}}^t+bG_{\mu}^{\mathcal{L}}({\boldsymbol{w}}_{t+1},\boldsymbol{p}_{t+1},\ell_{\mathcal{M}_1},f_{\mathcal{M}_2},\boldsymbol{u}_{[q]}),\\
\boldsymbol{z}_{\boldsymbol{p}}^{t+1}=&(1-b)\boldsymbol{z}_{\boldsymbol{p}}^t+bH(\boldsymbol{w}_{t+1},{\boldsymbol{p}}_{t+1},f_{\mathcal{M}_3}),
\end{align}
where $0<b<1$, $\boldsymbol{z}_{\boldsymbol{w}}^{1}=G_{\mu}^{\mathcal{L}}({\boldsymbol{w}}_{1},\boldsymbol{p}_{1},\ell_{\mathcal{M}_1},f_{\mathcal{M}_2},\boldsymbol{u}_{[q]})$ and $\boldsymbol{z}_{\boldsymbol{p}}^{1}=H(\boldsymbol{w}_{1},{\boldsymbol{p}}_{1},f_{\mathcal{M}_3})$. Here, $H(\boldsymbol{w}_{t+1},{\boldsymbol{p}}_{t+1},f_{\mathcal{M}_3})$ is calculated on the intermediate point $\boldsymbol{p}_{t+1}=(1-a)\boldsymbol{p}_{t}+a\hat{\boldsymbol{p}}_{t+1}$, where $0<a<1$ and $\hat{\boldsymbol{p}}_{t+1}$ is the distribution after updating and projecting onto the $\Delta^m$.

Then we use adaptive stepsizes to update $\boldsymbol{w}$ and $\boldsymbol{p}$. Specifically, we ensure the stepsizes are proportional to ${1}/({\sqrt{\|\boldsymbol{z}_{\boldsymbol{w}}^t\|_2}}+c)$ and ${1}/({\sqrt{\|\boldsymbol{z}_{\boldsymbol{p}}^t\|_2}}+c)$ \cite{liu2020adam,guo2021stochastic}, where $c>0$ is a small constant used to prevent the denominator from becoming $0$. Therefore, the update rules of $\boldsymbol{w}$ and $\boldsymbol{p}$ become
\begin{align}
\boldsymbol{w}_{t+1}&=\boldsymbol{w}_{t}-\eta_{\boldsymbol{w}}\dfrac{\boldsymbol{z}_{\boldsymbol{w}}^t}{\sqrt{\|\boldsymbol{z}_{\boldsymbol{w}}^t\|_2}+c},\\
\hat{\boldsymbol{p}}_{t+1}&=\mathcal{P}_{\Delta^m}(\boldsymbol{p}_t+\eta_{\boldsymbol{p}}\dfrac{\boldsymbol{z}_{\boldsymbol{p}}^t}{\sqrt{\|\boldsymbol{z}_{\boldsymbol{p}}^t\|_2}+c}).
\end{align}

These two key components of our method, \textit{i.e.,} extrapolation moving average and adaptive stepsize from the root norm of the momentum estimate, make our method enjoy two noticeable benefits: variance reduction of momentum estimate and adaptive convergence. The whole algorithm is presented in Algorithm \ref{alg:DSZOG3}. Since there exist two sources of randomness, we call our method doubly stochastic zeroth-order gradient method (DSZOG). Note that different from the algorithm in \cite{guo2021stochastic}, we use the adaptive step size method in both updating $\boldsymbol{w}$ and $\boldsymbol{p}$.  

\section{Convergence Analysis}
In this section, we discuss the convergence performance of our methods. The detailed proofs are given in our appendix.  

\subsection{Stationary point}
In this subsection, we first give the assumption about $\mathcal{L}$ and then give the definitions of the stationary point.
\begin{assumption}\label{assum:penalty_function}
	The objective function $\mathcal{L}(\boldsymbol{w},\boldsymbol{p})$ has the following properties:
	\begin{enumerate}
		\item $\mathcal{L}(\boldsymbol{w},\boldsymbol{p})$ is continuously differentiable in $\boldsymbol{w}$ and $\boldsymbol{p}$. $\mathcal{L}(\boldsymbol{w},\boldsymbol{p})$ is nonconvex with respect to $\boldsymbol{w}$, and $\mathcal{L}(\boldsymbol{w},\boldsymbol{p})$ is $\tau$-strongly concave with respect to $\boldsymbol{p}$.
		\item The function $g(\boldsymbol{w}):=\max_{\boldsymbol{p}}\mathcal{L}(\boldsymbol{w},\boldsymbol{p})$ is lower bounded, and $g(\boldsymbol{w})$ is $L_g$-Lipschitz continuous.
		
		\item When viewed as a function in $\mathbb{R}^{d+m}$, $\mathcal{L}(\boldsymbol{w},\boldsymbol{p})$ is $L$-gradient Lipschitz, ($L>0$), such that 
		$
		\|\nabla\mathcal{L}(\boldsymbol{w}_1,\boldsymbol{p}_1)-\nabla\mathcal{L}(\boldsymbol{w}_2,\boldsymbol{p}_2)\|_2
		\leq L\|(\boldsymbol{w}_1,\boldsymbol{p}_1)-(\boldsymbol{w}_2,\boldsymbol{p}_2)\|_2
		$.
	\end{enumerate}
\end{assumption}
This assumption is widely used in the convergence analysis of minimax problems \cite{wang2020zeroth,huang2020accelerated}. The first condition is used to detail the structure of $\mathcal{L}$ and the second condition is used to make the optimization problem well defined , and the third condition places a restriction on the degree of smoothness to be satisfied by the objective function.

Then, we discuss the definitions of stationary points and their relationships. For a general nonconvex constrained optimization problem, the stationary point \cite{lin2019inexact} is defined as follows,
\begin{definition}\label{def:kkt}
	$\boldsymbol{w}^*$ is said to be the stationary point of problem (\ref{constrained_problem}), if the following conditions holds,
	\begin{align}
	&\nabla_{\boldsymbol{w}} f_0(\boldsymbol{w}^*) + \sum_{j=1}^m\alpha_j^*\nabla_{\boldsymbol{w}}f_j(\boldsymbol{w}^*)=\boldsymbol{0},\\
	&f_j(\boldsymbol{w}^*)\leq 0, \\  
	&\alpha_j^*f_j(\boldsymbol{w}^*)=0, \quad \forall i\in\{1,\cdots,m\},
	\end{align}
	where $\boldsymbol{\alpha}^*=[\alpha_1,\cdots,\alpha_m]_t$ denotes the Lagrangian multiplier and $\alpha_j\geq0, \forall i=1,\cdots,m$.
\end{definition}
However, it is hard to compute a solution that satisfies the above conditions exactly \cite{lin2019inexact}. Therefore, finding the following $\epsilon$-stationary point \cite{lin2019inexact} is more practicable,
\begin{definition}\label{def:e-kkt}
	($\epsilon$-stationary) $\boldsymbol{w}^*$ is said to be the $\epsilon$-stationary point of problem (\ref{constrained_problem}), if there exists a vector $\boldsymbol{\alpha}^*\geq\boldsymbol{0}$, such that the following conditions hold,
	\begin{align}
	&\|\nabla_{\boldsymbol{w}} f_0(\boldsymbol{w}^*) + \sum_{j=1}^m\alpha_j^*\nabla_{\boldsymbol{w}}f_j(\boldsymbol{w}^*)\|_2^2\leq\epsilon_1^2, \\
	&\sum_{j=1}^m(\max\{f_j(\boldsymbol{w}^*),0\})^2\leq \epsilon_2^2,\\
	&\sum_{j=1}^m(\alpha_j f_j(\boldsymbol{w}^*))^2\leq \epsilon_3^2.
	\end{align}
\end{definition}

Since we reformulate the constrained problem as a minimax problem, here we give the definition of the approximation stationary point of the minimax problem and then show its relationship with Definition \ref{def:e-kkt}. According to \cite{wang2020zeroth}, we have the following definition,
\begin{definition}\label{def:minimax}
	A point $(\boldsymbol{w}^*,\boldsymbol{p}^*)$ is called the $\epsilon$-stationary point of problem $\min_{\boldsymbol{w}}\max_{\boldsymbol{p}\in\Delta^m}\mathcal{L}(\boldsymbol{w},\boldsymbol{p})$ if it satisfies the conditions: $\|\nabla_{\boldsymbol{w}}\mathcal{L}(\boldsymbol{w},\boldsymbol{p})\|_2^2\leq\epsilon^2$ and $\|\nabla_{\boldsymbol{p}}\mathcal{L}(\boldsymbol{w},\boldsymbol{p})\|_2^2\leq\epsilon^2$.
\end{definition}
In addition, we have the following Proposition between definition \ref{def:e-kkt} and definition \ref{def:minimax}.
\begin{proposition}\label{prop1}
	If Assumption \ref{assum:penalty_function} holds, $\sqrt{\dfrac{2m\epsilon^2+2m^2\lambda^2}{\beta^2}}\leq\epsilon_2^2$ and $(\boldsymbol{w}^*,\boldsymbol{p}^*)$ is the $\epsilon$-stationary point defined in Definition \ref{def:minimax} of the problem $\min_{\boldsymbol{w}}\max_{\boldsymbol{p}\in\Delta^m}\mathcal{L}(\boldsymbol{w},\boldsymbol{\alpha})$, then $\boldsymbol{w}^*$ is the $\epsilon$-stationary point defined in Definition \ref{def:e-kkt} of the constrained problem \ref{constrained_problem}.
\end{proposition}

As proposed in \cite{wang2020zeroth}, the minimax problem \ref{mini_opt} is equivalent to the following minimization problem:
\begin{align}\label{min_problem}
\min_{\boldsymbol{w}}\left\{ g(\boldsymbol{w}):=\max_{\boldsymbol{p}\in\Delta^m}\mathcal{L}(\boldsymbol{w},\boldsymbol{p})=\mathcal{L}(\boldsymbol{w},\boldsymbol{p}^*(\boldsymbol{w})) \right\},
\end{align}
where $\boldsymbol{p}^*(\boldsymbol{w})=\arg\max_{\boldsymbol{p}}\mathcal{L}(\boldsymbol{w},\boldsymbol{p})$.
Here, we give stationary point the minimization problem \ref{min_problem}  and its relationship with Definition \ref{def:minimax} as follows,
\begin{definition}\label{def:min_stationary}
	We call $\boldsymbol{w}^*$ an $\epsilon$-stationary point of a differentiable function $g(\boldsymbol{w})$, if $\|\nabla g(\boldsymbol{w}^*)\|_2\leq \epsilon$.
\end{definition}
\begin{proposition}\label{prop2}
	Under Assumption \ref{assum:penalty_function}, if a point $\boldsymbol{w}'$ is an $\epsilon$-stationary point in terms of Definition \ref{def:min_stationary}, then an $\epsilon$-stationary point $\boldsymbol{w}^*,\boldsymbol{p}^*$ in terms of Definition \ref{def:minimax} can be obtained. 
\end{proposition}

\begin{remark}
	According to Proposition \ref{prop1} and Proposition \ref{prop2}, we have that once we find the $\epsilon$-stationary point in terms of Definition \ref{def:min_stationary}, then we can get the $\epsilon$-stationary point in terms of Definition \ref{def:e-kkt}.
\end{remark}

\subsection{Convergence Rate of the Accelerated Method}
In this subsection, we discuss the convergence performance of our algorithms. Here, we give several assumptions used in our analysis.
\begin{assumption}\label{assum:bound_step}
	We have $c_{1,l}\leq\dfrac{1}{\sqrt{\|\boldsymbol{z}_{\boldsymbol{w}}^t\|_2}+c}\leq c_{1,u}$ and $c_{2,l}\leq\dfrac{1}{\sqrt{\|\boldsymbol{z}_{\boldsymbol{p}}^t\|_2}+c}\leq c_{2,u}$, where $c$ is a constant.
\end{assumption}
This assumption is used to bound the step size scaling factor which is widely used in \cite{huang2021super,guo2021stochastic,huang2021biadam}.
\begin{assumption}\label{assum:mean_variance}
	For any $\boldsymbol{w}\in\mathbb{R}^d$, the following properties holds,
	\begin{align}
	&\mathbb{E}[G_{\mu}^{\mathcal{L}}(\boldsymbol{w},\boldsymbol{p},\ell_{\mathcal{M}_1},f_{\mathcal{M}_2},\boldsymbol{u}_{[q]})]=\nabla_{\boldsymbol{w}}\mathcal{L}(\boldsymbol{w},\boldsymbol{p}),\nonumber\\
	&\mathbb{E}[H(\boldsymbol{w},\boldsymbol{p},f_{\mathcal{M}_3})]=\nabla_{\boldsymbol{p}}\mathcal{L}(\boldsymbol{w},\boldsymbol{p}),\nonumber\\
	&\mathbb{E}[\|G_{\mu}^{\mathcal{L}}(\boldsymbol{w}_1,\boldsymbol{p}_1,\ell_{\mathcal{M}_1},f_{\mathcal{M}_2},\boldsymbol{u}_{[q]})-\nabla_{\boldsymbol{w}}\mathcal{L}(\boldsymbol{w}_1,\boldsymbol{p}_1)\|_2]\leq\sigma_1^2,\nonumber\\ &\mathbb{E}[\|H(\boldsymbol{w},\boldsymbol{p},f_{\mathcal{M}_3})-\nabla_{\boldsymbol{w}}\mathcal{L}(\boldsymbol{w}_t,\boldsymbol{p}_t)\|_2]\leq\sigma_2^2.\nonumber
	\end{align}.
\end{assumption}
This assumptions is used to show our estimator is the unbiased estimation. Based on above assumptions, we can derive the following lemmas which are useful in our convergence analysis.
\begin{lemma}\label{lemma:bound_of_g_in_ASZOG}
 (\textbf{Descent in the function value.})	Under Assumptions \ref{assum:penalty_function} and \ref{assum:bound_step}, if $\eta_{w}L\leq\dfrac{c_{1,l}}{2c_{1,u}^2}$, we have
	\begin{align}
	&g(\boldsymbol{w}_{t+1})\nonumber\\
	\leq&g(\boldsymbol{w}_t)-\dfrac{\eta_wc_{1,l}}{2}\|\nabla g(\boldsymbol{w}_t)\|_2^2-\dfrac{\eta_wc_{1,l}}{4}\|\boldsymbol{z}_{\boldsymbol{w}}^t\|_2^2\nonumber\\
	&+\eta_wc_{1,l}\dfrac{\mu^2 L^2(d+3)^3}{2}+\eta_wc_{1,l}L^2\|\boldsymbol{p}^*(\boldsymbol{w}_t)-\boldsymbol{p}_t\|_2^2\nonumber\\
	&+\eta_wc_{1,l}\|\nabla_{\boldsymbol{w}}\mathcal{L}(\boldsymbol{w}_t,\boldsymbol{p}_t)-\boldsymbol{z}_{\boldsymbol{w}}^t\|_2^2.\nonumber
	\end{align}
\end{lemma}
\begin{lemma}\label{lemma:bound_of_p_in_ASZOGD}
	(\textbf{Descent in the iterates of the probability.}) Under Assumptions \ref{assum:penalty_function}, \ref{assum:bound_step} and \ref{assum:mean_variance}, if $a\leq1$ and $\eta_{p}\leq\dfrac{1}{3c_{2.l}L}$, we have
	\begin{align}
	&\|\boldsymbol{p}_{t+1}-\boldsymbol{p}^*(\boldsymbol{w}_{t+1})\|_2^2\nonumber\\
	\leq&-\dfrac{1}{4a}\|\boldsymbol{p}_t-{\boldsymbol{p}}_{t+1}\|_2^2+\dfrac{8a\eta_pc_{2,l}}{\tau}\|\nabla_{\boldsymbol{p}}\mathcal{L}(\boldsymbol{w}_t,\boldsymbol{p}_t)-\boldsymbol{z}_{\boldsymbol{p}}^t\|_2^2
	\nonumber\\
	&+(1-\dfrac{\tau a \eta_{p}c_{2,l}}{4})\|\boldsymbol{p}^*(\boldsymbol{w}_{t})-\boldsymbol{p}_t\|_2^2\nonumber\\
	&+\dfrac{8L^2_g}{\tau a \eta_{p}c_{2,l}}\|\boldsymbol{w}_{t}-\boldsymbol{w}_{t+1}\|_2^2.\nonumber
	\end{align}
\end{lemma}
\begin{lemma}\label{lemma:bound_of_z_in_ASZOGD}
	(\textbf{Descent in the gradient estimation error.}) Under Assumptions \ref{assum:penalty_function}, \ref{assum:bound_step} and \ref{assum:mean_variance}, if $b\in(0,1)$, we have
	\begin{align}
	&\mathbb{E}[\| \nabla_{\boldsymbol{p}}\mathcal{L}(\boldsymbol{w}_{t+1},\boldsymbol{p}_{t+1}) -\boldsymbol{z}_{\boldsymbol{p}}^{t+1} \|_2^2]\nonumber\\
	\leq&(1-b)\mathbb{E}[\|\nabla_{\boldsymbol{p}}\mathcal{L}(\boldsymbol{w}_t,\boldsymbol{p}_t)-\boldsymbol{z}_{\boldsymbol{p}}^t\|_2^2]+\dfrac{1}{b}L^2[\|\boldsymbol{w}_{t+1}-\boldsymbol{w}_t\|_2^2\nonumber\\
	&+\|\boldsymbol{p}_{t+1}-\boldsymbol{p}_t\|_2^2]+b^2\sigma_2^2,\nonumber\\
	&\mathbb{E}[\| \nabla_{\boldsymbol{w}}\mathcal{L}(\boldsymbol{w}_{t+1},\boldsymbol{p}_{t+1}) -\boldsymbol{z}_{\boldsymbol{w}}^{t+1} \|_2^2]\nonumber\\
	\leq&(1-b)\mathbb{E}[\|\nabla_{\boldsymbol{w}}\mathcal{L}(\boldsymbol{w}_t,\boldsymbol{p}_t)-\boldsymbol{z}_{\boldsymbol{w}}^t\|_2^2]+\dfrac{1}{b}L^2[\|\boldsymbol{w}_{t+1}-\boldsymbol{w}_t\|_2^2\nonumber\\
	&+\|\boldsymbol{p}_{t+1}-\boldsymbol{p}_t\|_2^2]+b^2\sigma_1^2.\nonumber
	\end{align}
\end{lemma}
Then, following the framework in \cite{guo2021stochastic,wang2018stochastic,huang2020accelerated} and utilizing the above lemmas, we have the following theorem,
\begin{theorem}\label{theorem:convergence_rate_of_ASZOGD}
	Under Assumptions \ref{assum:penalty_function}, \ref{assum:bound_step} and \ref{assum:mean_variance}, if $a\in(0,1]$, ,$\boldsymbol{p}^*(w_{1})=\boldsymbol{p}_1$, $\boldsymbol{z}_{\boldsymbol{p}}^1=H(\boldsymbol{w}_{t},\boldsymbol{p}_{t},f_{\mathcal{M}_3})$, $\boldsymbol{z}_{\boldsymbol{w}}^1=G_{\mu}^{\mathcal{L}}(\boldsymbol{w}_{t},\boldsymbol{p}_{t},\ell_{\mathcal{M}_1},f_{\mathcal{M}_2},\boldsymbol{u}_{[q]})$, $0<\eta_{p}\leq\min\{\dfrac{1}{3c_{2.l}L},\dfrac{b^2}{\tau a^2c_{2,l}},\dfrac{\tau b^2}{32L^2a^2c_{2,l}},1\}$, $0<\eta_{w}^2\leq\min\{\dfrac{c_{1,l}^2}{4Lc_{1,u}^4},\dfrac{b^2}{4c_{1,u}^2L^2},\dfrac{\tau^2a^2\eta_{p}^2c_{2,l}^2}{128L_g^2L^2c_{1,u}},\dfrac{\tau^2b^2}{128L^4c_{1,u}^2},1\}$, $\mu\leq\dfrac{\epsilon}{L(d+3)^{3/2}}$, $0<b\leq\min\{\dfrac{\epsilon^2}{2\sigma_1^2},\dfrac{\tau^2\epsilon^2}{64\sigma_2^2L^2},1\}$ and $T\geq\max\{\dfrac{2(g(\boldsymbol{w}_1)-g(\boldsymbol{w}_T))}{\epsilon^2\eta_{w}c_{1,l}},\dfrac{2\sigma_1^2}{\epsilon^2b},\dfrac{64\sigma_2^2L^2}{\epsilon^2\tau^2 b}\}$, we have
	\begin{align}
	\dfrac{1}{T}\mathbb{E}[\sum_{t=1}^{T}\|\nabla g(\boldsymbol{w}_t)\|_2^2]
	\leq\epsilon^2.
	\end{align}
\end{theorem}
\begin{remark}
	By choosing $b=\mathcal{O}({\epsilon^2}/{\kappa^2})$, $\eta_{\boldsymbol{w}}=\mathcal{O}({\epsilon^2}/{\kappa^6})$, $\eta_{\boldsymbol{p}}$=$\mathcal{O}({\epsilon^2}/{\kappa^4})$ and $T=\mathcal{O}({\kappa^6}/{\epsilon^4})$, our proposed DSZOG can converge to the $\epsilon$-stationary point defined in Definition \ref{def:min_stationary}. Then, based on Proposition \ref{prop1} and Proposition \ref{prop2}, we can derive that our method can converge to the $\epsilon$-stationary point of the original constrained problem (\ref{constrained_problem}) defined in Definition \ref{def:e-kkt}.
\end{remark}

\begin{table*}[]
	\centering
		\setlength{\tabcolsep}{5.2mm}
	\caption{Test accuracy (\%) of all the methods in classification with pairwise constraints.}
	\label{tab:classification_pairwise}
	\begin{tabular}{llllll}
		\toprule
		Data& DSZOG                & ZOSCGD        & ZOPSGD &AccZOMDA& AccSZOFW                   \\
		\hline
		a9a    & $\textbf{75.90}\pm0.26$   & $75.35\pm0.13$ & $75.37\pm0.19$ &$75.52\pm0.21$&$75.22\pm0.12$ \\
		\hline
		w8a & $\textbf{89.94}\pm0.28$  & $83.53\pm0.58$ & $89.02\pm0.97$&$89.14\pm0.23$&$89.34\pm0.33$  \\
		\hline
		gen & $\textbf{82.33}\pm0.76$  & $66.33\pm0.07$ & $66.83\pm0.57$&$72.84\pm0.45$&$72.03\pm0.23$  \\
		\hline
		svm   & $\textbf{79.56}\pm0.49$  & $71.21\pm0.57$ & $78.63\pm0.26$&$79.01\pm0.21$&$71.88\pm0.34$ \\
		\bottomrule
	\end{tabular}	
\end{table*}
\begin{table*}[]
	\centering
		\setlength{\tabcolsep}{5mm}
	\caption{Test accuracy (\%) of all the methods in classification with fairness constraints.}
	\label{tab:fairness}
	\begin{tabular}{llllll}
		\toprule
		Data& DSZOG                & ZOSCGD        & ZOPSGD    &AccZOMDA& AccSZOFW               \\
		\hline
		D1    & $\textbf{87.33}\pm0.38$   & $51.08\pm0.57$ & $59.16\pm0.37$  &$66.33\pm0.19$&$55.23\pm0.46$\\
		\hline
		D2 & $\textbf{84.75}\pm0.25$  & $69.70\pm0.24$ & $68.00\pm0.54$ &$69.55\pm0.29$&$76.13\pm0.45$ \\
		\hline
		D3 & $\textbf{83.58}\pm0.14$ & $66.33\pm0.30$ & $66.84\pm0.57$&$66.35\pm0.45$&$60.47\pm0.66$  \\
		\hline
		D4   & $\textbf{64.91}\pm0.94$  & $52.16\pm0.38$ & $55.25\pm0.90$&$55.40\pm0.51$&$54.86\pm0.43$ \\
		\bottomrule
	\end{tabular}

\end{table*}

\section{Experiments}

\subsection{Experimental Setup}
In this subsection, we summarized the baselines used in our experiments as follows,
\begin{enumerate}
	\item \textbf{ZOPSGD}\cite{liu2018zeroth_pro}. In each iteration, ZOPSGD calculates the stochastic zeroth-order gradient of $f_0$ to update the parameters and then solves a constrained quadratic problem to project the solution into the feasible set. 

	\item \textbf{ZOSCGD}\cite{balasubramanian2018zeroth}. In each iteration, ZOSCGD calculates the stochastic zeroth-order gradient of $f_0$ and then uses the conditional gradient method to update the parameters by solving a constrained linear problem.
	\item \textbf{AccZOMDA}\cite{huang2020accelerated}. In each iteration, AccZOMDA uses the momentum-based variance reduce technique of STORM \cite{cutkosky2019momentum} to estimate the stochastic zeroth-order gradients, and then solves a constrained quadratic problem to project the solution into the feasible set.
	\item \textbf{AccSZOFW}\cite{huang2020accelerated2}.In each iteration, AccSZOFW uses the variance reduced technique of SPIDER \cite{fang2018spider} to calculate the stochastic zeroth-order gradient of $f_0$ and then uses the conditional gradient method to update the parameters by solving a constrained linear problem.
\end{enumerate}

\subsection{Applications}
\begin{table}
	\setlength{\tabcolsep}{7mm}
	\centering
	\caption{Datasets used in classification with pairwise constraints (We give the approximate size of constraints).}
	\label{tab:data1}
	\begin{tabular}{lcc}
		\toprule
		Data      & Dimension & Constriants \\\hline
		w8a       & 300       & $\simeq$8000        \\
		a9a       & 123       & $\simeq$40000       \\
		gen       & 50        & $\simeq$60000       \\
		svm & 22        & $\simeq$40000      \\
		\bottomrule
	\end{tabular}
\end{table}

\begin{figure*}[]
	\centering
	
	\begin{subfigure}[b]{0.245\textwidth}
		\centering
		\includegraphics[width=1.8in]{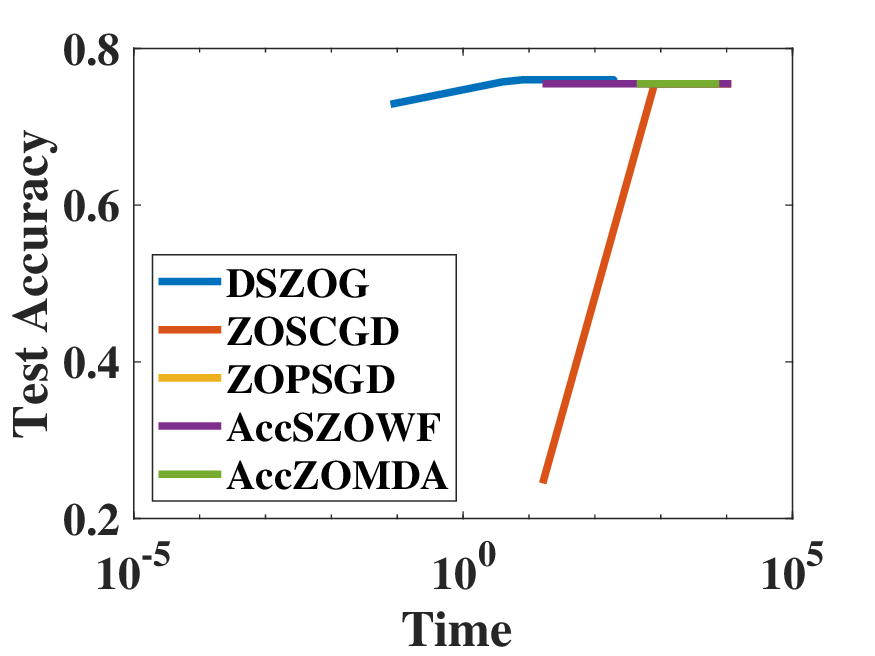}
		\caption{a9a}
	\end{subfigure}
	\begin{subfigure}[b]{0.245\textwidth}
		\centering
		\includegraphics[width=1.8in]{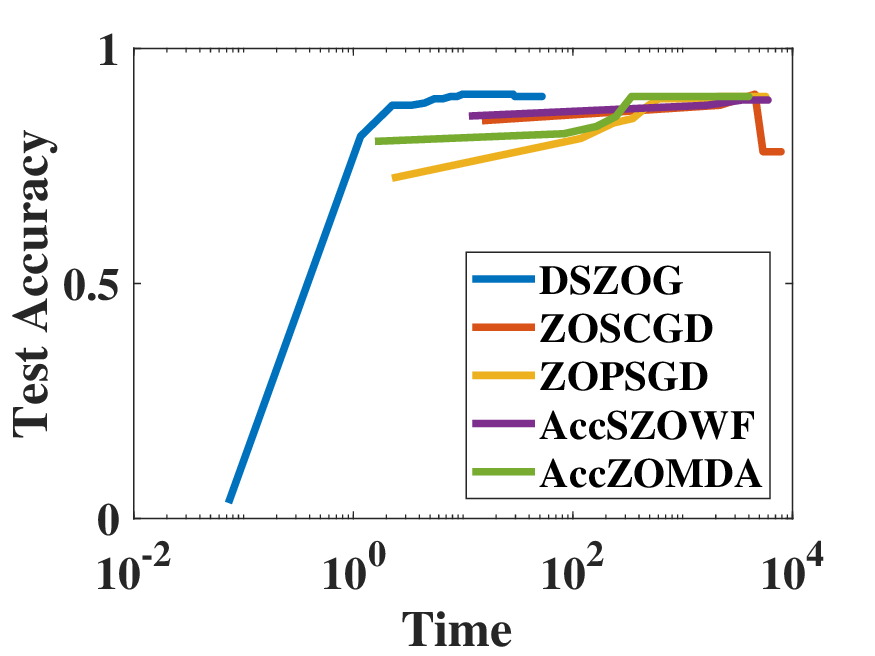}
		\caption{w8a}
	\end{subfigure}
	\begin{subfigure}[b]{0.245\textwidth}
		\centering
		\includegraphics[width=1.8in]{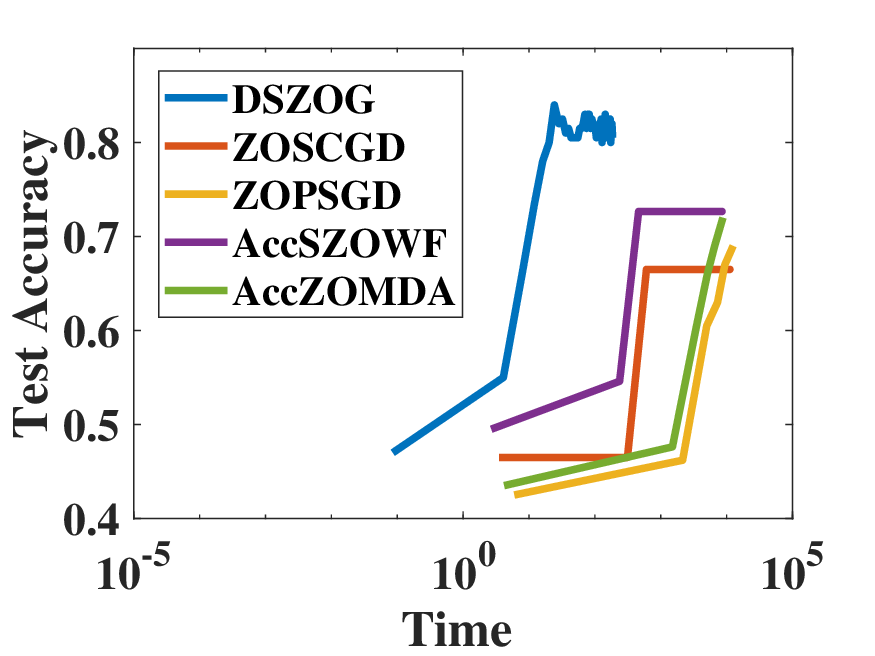}
		\caption{gen}
	\end{subfigure}	
	\begin{subfigure}[b]{0.245\textwidth}
		\centering
		\includegraphics[width=1.8in]{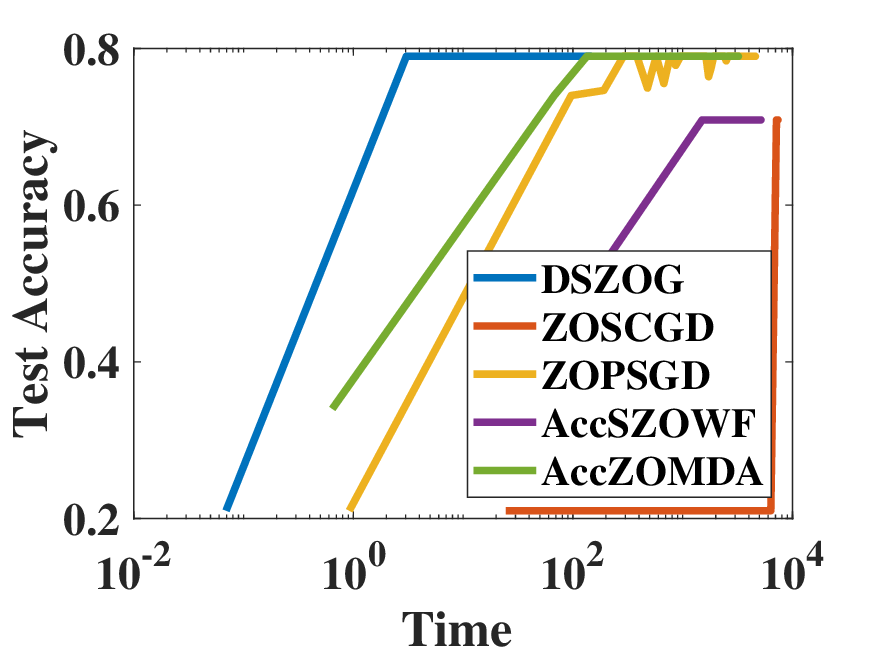}
		\caption{svm}
	\end{subfigure}
	
	\caption{Test accuracy against training time of all the methods in classification with pairwise constraints (We stop the algorithms if the training time is more than 10000 seconds).}
	\label{fig:classification_pairwise}
\end{figure*}

\begin{figure*}[]
	\centering
	\small
	\begin{subfigure}[b]{0.245\textwidth}
		\centering
		\includegraphics[width=1.8in]{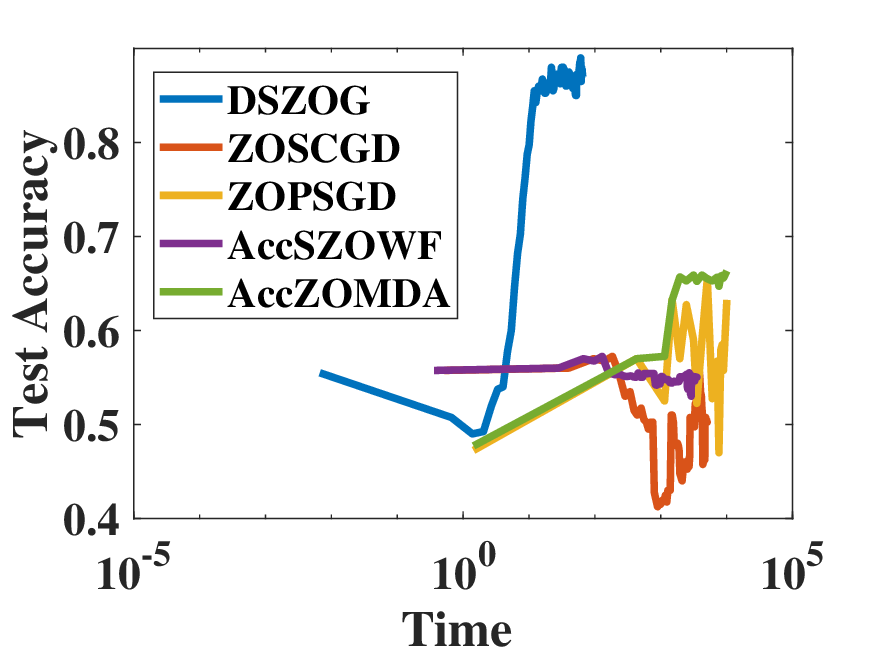}
		\caption{D1}
	\end{subfigure}
	\begin{subfigure}[b]{0.245\textwidth}
		\centering
		\includegraphics[width=1.8in]{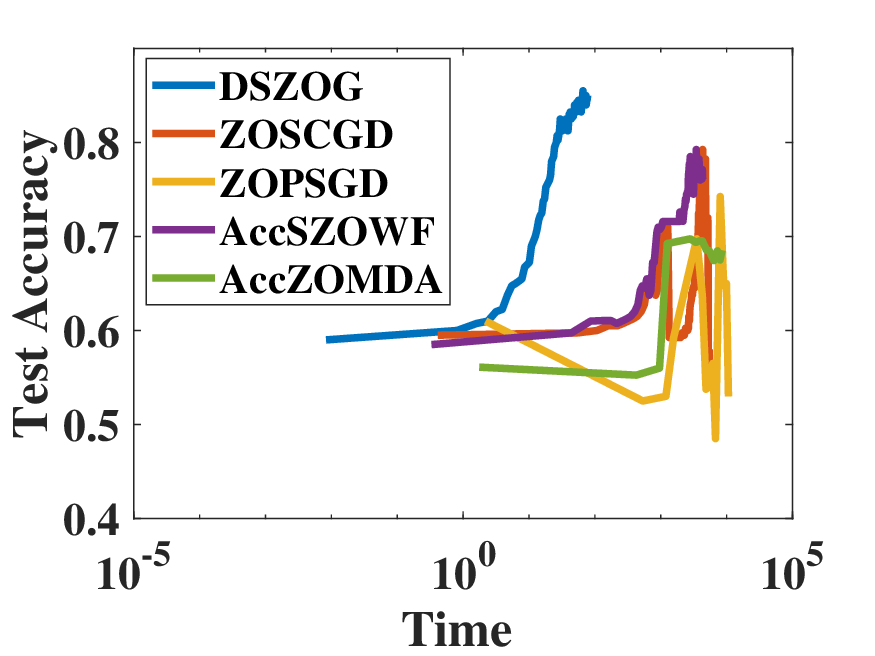}
		\caption{D2}
	\end{subfigure}
	\begin{subfigure}[b]{0.245\textwidth}
		\centering
		\includegraphics[width=1.8in]{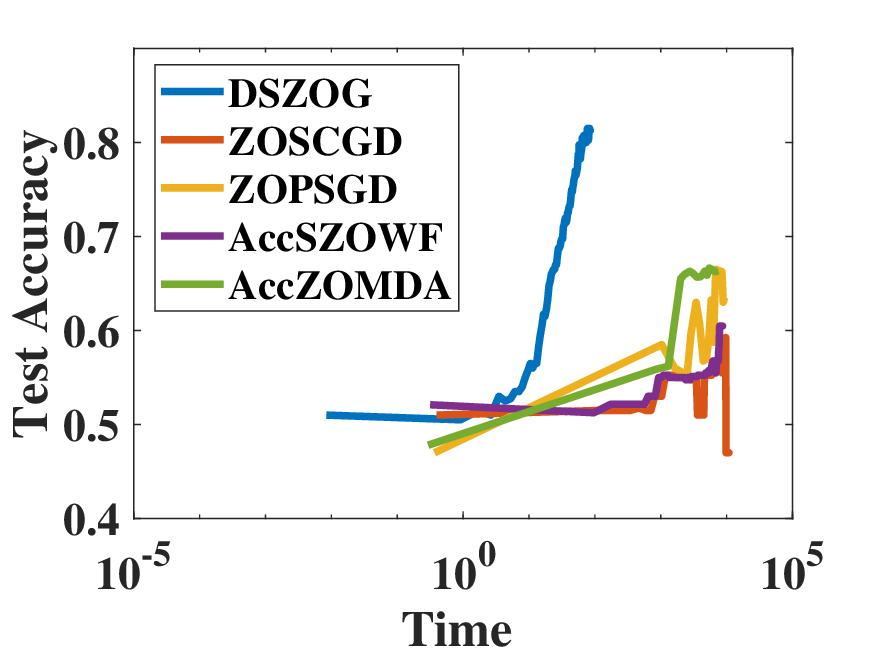}
		\caption{D3}
	\end{subfigure}	
	\begin{subfigure}[b]{0.245\textwidth}
		\centering
		\includegraphics[width=1.8in]{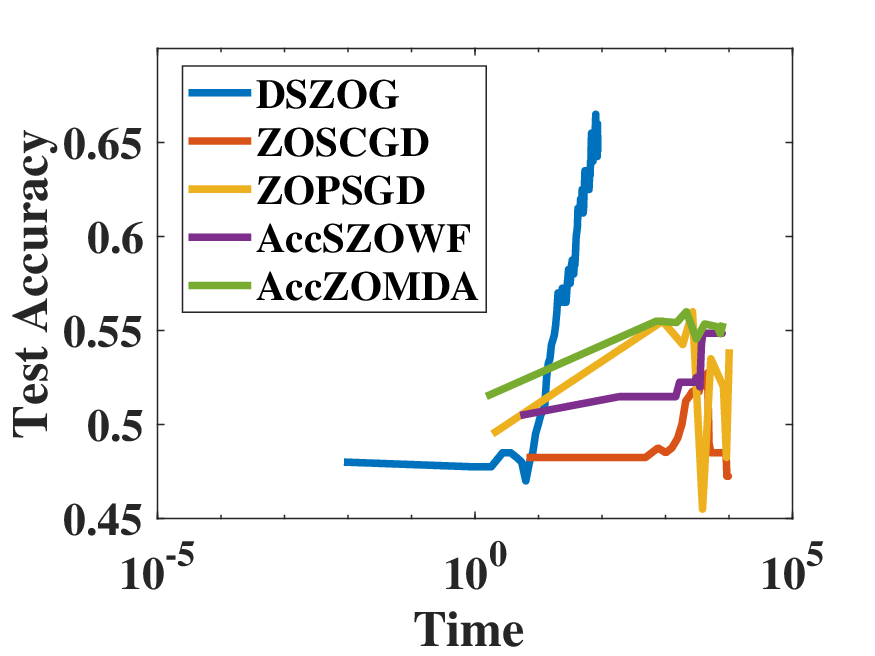}
		\caption{D4}
	\end{subfigure}
	
	\caption{Test accuracy against training time of all the methods in classification with fairness constraints (We stop the algorithms if the training time is more than 10000 seconds).}
	\label{fig:fairness}
\end{figure*}

In this subsection, we give the introduction of the applications used in our experiments.

\noindent \textbf{Classification with Pairwise Constraints}
We evaluate the performance of all the methods on the binary classification with pairwise constraints learning problem.  Given a set of training samples $\mathcal{D}=\{(\boldsymbol{x}_i,y_i)\}_{i=1}^n$, where $\boldsymbol{x}_i\in\mathbb{R}^d$ and $y_i\in\{+1,-1\}$. In this task, we learn a linear model $h(\boldsymbol{x},\boldsymbol{w})=\boldsymbol{x}^T\boldsymbol{w}$ to classify the dataset and ensure the any positive sample $\boldsymbol{x}_i^{+}\in\mathcal{D}^{+}:=\{(\boldsymbol{x}_i,+1)\}_{i=1}^{n_p}$ has larger function value than the negative sample $\boldsymbol{x}_j^{-}\in\mathcal{D}^{-}:=\{(\boldsymbol{x}_j,+1)\}_{i=1}^{n_n}$, where $n_p$ and $n_n$ denotes the number of positive samples and negative samples, respectively. Then, we can formulate this problem as follows,
\begin{align}
\min_{\boldsymbol{w}}&\ \dfrac{1}{n}\sum_{i=1}^n\ell(h(\boldsymbol{x}_i,\boldsymbol{w}),y_i),\\
s.t. & \ h(\boldsymbol{x}_i^{+},\boldsymbol{w})-h(\boldsymbol{x}_j^{-},\boldsymbol{w})\geq 0,\nonumber\\
& \ \forall\boldsymbol{x}_i^{+}\in\mathcal{D}^{+}\ \boldsymbol{x}_j^{-}\in\mathcal{D}^{-},\nonumber
\end{align}
where $\ell(u,v) = c^2(1-\exp(-\dfrac{(v-u)^2}{c^2})) $ is viewed as a black-box function. We summarized the datasets used in this application in Table \ref{tab:data1}. We randomly sample $1000$ data samples from the original datasets, and then divide all the datasets into 3 parts, i.e., $50\%$ for training, $30\%$ for testing and $20\%$ for validation. We fix the batch size of data sample at $128$ for all the methods and $|\mathcal{M}_2|=|\mathcal{M}_3|=128$. The learning rates of all the methods are chosen from $\{0.01,0.001,0.0001\}$. In our methods, the penalty parameter $\beta$ is chosen from $\{0.1,1,10\}$, $a$ and $b$ are chosen from $\{0.1,0.5,0.9\}$ on the validation sets.

\noindent\textbf{Classification with Fairness Constraints.}
In this problem, we consider the binary classification problem with a large amount of fairness constraints \cite{zafar2017fairness}. Given a set of training samples $\mathcal{D}=\{(\boldsymbol{x}_i,y_i)\}_{i=1}^n$, where $\boldsymbol{x}_i\in\mathbb{R}^d$ and $y_i\in\{-1,+1\}$. In this task, we learn a linear model $h(\boldsymbol{x},\boldsymbol{w})=\boldsymbol{x}^T\boldsymbol{w}$. Assume that each sample has an associate sensitive feature vector $\boldsymbol{z}\in\mathbb{R}^r$. We denote $z_{ij}\in\{0,1\}$ as the $j$-th sensitive feature of $i$-th sample. The classifier $h$ cannot use the protected characteristic $\boldsymbol{z}$ at decision time, as it will constitute an unfair treatment. A number of metrics have been used to determine how fair a classifier is with respect to the sensitive features. According to \cite{zafar2017fairness}, the fair classification problems can be formulated as follows,
\begin{align}
\min_{\boldsymbol{w}} &\dfrac{1}{n}\sum_{i=1}^n\ell(h(\boldsymbol{x}_i,\boldsymbol{w}),y_i),\\
s.t. &\ \dfrac{1}{n}\sum_{i=1}^n(z_{ij}-\bar{z}_j)g(y_i,\boldsymbol{x}_i)\leq c,\nonumber\\
 &\ 
\dfrac{1}{n}\sum_{i=1}^n(z_{ij}-\bar{z}_j)g(y_i,\boldsymbol{x}_i)\geq -c, \nonumber
\end{align}
where $j=1,\cdots,r$, $\ell(u,v)$ denotes the loss functions, $c$ is the covariance threshold which specifies an upper bound on the covariance between the sensitive attributes $\boldsymbol{z}$ and the signed distance $g(y,\boldsymbol{x})$. We use the hinge loss $\ell(u,v)=\max\{1-uv,0\}$ in this experiment and we view it as a black-box function. In addition, we use the following two functions to build the fairness constraints,
$
g(y,\boldsymbol{x})=\left\{\begin{matrix}
&\min\{0,\dfrac{1+y}{2}yh(\boldsymbol{x},\boldsymbol{w})\}\\
&\min\{0,\dfrac{1-y}{2}h(\boldsymbol{x},\boldsymbol{w})\}\\
\end{matrix}\right.
$. Since the datasets with multiple sensitive features are difficult to find,  we generate 4 datasets with $2000$ samples in this task and summarize them in Table \ref{tab:data2}. For each dataset, we randomly choose several features to be the sensitive features, and then separate them into 3 parts, i.e., $50\%$ for training, $30\%$ for testing and $20\%$ for validation. We fix the batch size of data sample at $128$ for all the methods and $|\mathcal{M}_2|=|\mathcal{M}_3|=10$. The learning rates of all the methods are chosen from $\{0.01,0.001,0.0001\}$. For our methods, the penalty parameter $\beta$ is chosen from $\{0.1,1,10\}$, $a$ and $b$ are chosen from $\{0.1,0.5,0.9\}$ on the validation sets.
\begin{table}
	\centering
	\caption{Datasets used in classification with fairness constraints.}
	\label{tab:data2}
	\begin{tabular}{lccc}
		\toprule
		Data      & Dimension & Sensitive Features& Constraints \\\hline
		D1       & 100       &10& 40       \\
		D2       & 200       &20&  80      \\
		D3       & 300        &20& 80       \\
		D4  	& 400        &20&  80     \\
		\bottomrule
	\end{tabular}
\end{table}
We run all the methods 10 times on a 3990x workstation.

\subsection{Results and Discussion}

We present the results in Figures \ref{fig:classification_pairwise}, \ref{fig:fairness} and Tables \ref{tab:classification_pairwise}, \ref{tab:fairness}. Note that for ZOSCGD, ZOPSGD, AccSZOFW and AccZOMDA, if the training time is larger than 10000 seconds, the algorithms are stopped. From Tables \ref{tab:classification_pairwise} and \ref{tab:fairness}, we can find that our methods DSZOG has the highest test accuracy in most cases in both two applications. In addition, from Figures \ref{fig:classification_pairwise} and \ref{fig:fairness}, we can find that our methods are faster than other methods. This is because all the other methods need to solve a subproblem with a large number of constraints in each iteration and the existing Python package cannot efficiently deal with such a problem. What's worse,  ZOSCGD, ZOPSGD, AccSZOFW and AccZOMDA focus on solving the problem with convex constraints while the constraints in the fairness problem are nonconvex. This makes ZOSCGD, ZOPSGD, AccSZOFW and AccZOMDA cannot find the stationary point. However, by using the penalty framework, our methods can still converge to the stationary point when the constraints are nonconvex. In addition, by using a stochastic manner on the constraint, our method can efficiently deal with a large number of constraints. 
All these results demonstrate that our method is superior to ZOSCGD and ZOPSGD in the heavily constrained nonconvex problem.

\section{Conclusion}
In this paper, we propose two efficient ZO method to solve the heavily constrained nonconvex black-box problem, i.e., DSZOG. We add an additional stochastic layer into the constraint to estimate the zeroth-order gradients. In addition, momentum and adaptive step size is also used in our method. We give the convergence analysis of our proposed method. The experimental results on two applications demonstrate the superiority of our method in terms of accuracy and training time .

\section*{Acknowledgments}
This work was partially supported by the National Natural Science Foundation of China under Grant 62076138, Postgraduate Research $\&$ Practice Innovation Program of Jiangsu Province under Grant KYCX21$\_$0999.

\bibliography{example_paper}

\begin{thebibliography}{47}
\providecommand{\natexlab}[1]{#1}
\providecommand{\url}[1]{\texttt{#1}}
\expandafter\ifx\csname urlstyle\endcsname\relax
  \providecommand{\doi}[1]{doi: #1}\else
  \providecommand{\doi}{doi: \begingroup \urlstyle{rm}\Url}\fi

\bibitem[Agarwal et~al.(2010)Agarwal, Dekel, and Xiao]{agarwal2010optimal}
Agarwal, A., Dekel, O., and Xiao, L.
\newblock Optimal algorithms for online convex optimization with multi-point
  bandit feedback.
\newblock In \emph{COLT}, pp.\  28--40. Citeseer, 2010.

\bibitem[Allen-Zhu(2017)]{allen2017katyusha}
Allen-Zhu, Z.
\newblock Katyusha: The first direct acceleration of stochastic gradient
  methods.
\newblock \emph{The Journal of Machine Learning Research}, 18\penalty0
  (1):\penalty0 8194--8244, 2017.

\bibitem[Balasubramanian \& Ghadimi(2018)Balasubramanian and
  Ghadimi]{balasubramanian2018zeroth}
Balasubramanian, K. and Ghadimi, S.
\newblock Zeroth-order (non)-convex stochastic optimization via conditional
  gradient and gradient updates.
\newblock In \emph{Proceedings of the 32nd International Conference on Neural
  Information Processing Systems}, pp.\  3459--3468, 2018.

\bibitem[Chen et~al.(2020)Chen, Zhou, Yi, and Gu]{chen2020frank}
Chen, J., Zhou, D., Yi, J., and Gu, Q.
\newblock A frank-wolfe framework for efficient and effective adversarial
  attacks.
\newblock In \emph{Proceedings of the AAAI Conference on Artificial
  Intelligence}, volume~34, pp.\  3486--3494, 2020.

\bibitem[Chen et~al.(2017)Chen, Zhang, Sharma, Yi, and Hsieh]{chen2017zoo}
Chen, P.-Y., Zhang, H., Sharma, Y., Yi, J., and Hsieh, C.-J.
\newblock Zoo: Zeroth order optimization based black-box attacks to deep neural
  networks without training substitute models.
\newblock In \emph{Proceedings of the 10th ACM workshop on artificial
  intelligence and security}, pp.\  15--26, 2017.

\bibitem[Choromanski et~al.(2018)Choromanski, Rowland, Sindhwani, Turner, and
  Weller]{choromanski2018structured}
Choromanski, K., Rowland, M., Sindhwani, V., Turner, R., and Weller, A.
\newblock Structured evolution with compact architectures for scalable policy
  optimization.
\newblock In \emph{International Conference on Machine Learning}, pp.\
  970--978. PMLR, 2018.

\bibitem[Clarkson et~al.(2012)Clarkson, Hazan, and
  Woodruff]{clarkson2012sublinear}
Clarkson, K.~L., Hazan, E., and Woodruff, D.~P.
\newblock Sublinear optimization for machine learning.
\newblock \emph{Journal of the ACM (JACM)}, 59\penalty0 (5):\penalty0 1--49,
  2012.

\bibitem[Cotter et~al.(2016)Cotter, Gupta, and Pfeifer]{cotter2016light}
Cotter, A., Gupta, M., and Pfeifer, J.
\newblock A light touch for heavily constrained sgd.
\newblock In \emph{Conference on Learning Theory}, pp.\  729--771. PMLR, 2016.

\bibitem[Cutkosky \& Mehta(2020)Cutkosky and Mehta]{cutkosky2020momentum}
Cutkosky, A. and Mehta, H.
\newblock Momentum improves normalized sgd.
\newblock In \emph{International Conference on Machine Learning}, pp.\
  2260--2268. PMLR, 2020.

\bibitem[Cutkosky \& Orabona(2019)Cutkosky and Orabona]{cutkosky2019momentum}
Cutkosky, A. and Orabona, F.
\newblock Momentum-based variance reduction in non-convex sgd.
\newblock \emph{Advances in Neural Information Processing Systems},
  32:\penalty0 15236--15245, 2019.

\bibitem[Defazio et~al.(2014)Defazio, Bach, and
  Lacoste-Julien]{defazio2014saga}
Defazio, A., Bach, F., and Lacoste-Julien, S.
\newblock Saga: A fast incremental gradient method with support for
  non-strongly convex composite objectives.
\newblock In \emph{Advances in neural information processing systems}, pp.\
  1646--1654, 2014.

\bibitem[Duchi et~al.(2015)Duchi, Jordan, Wainwright, and
  Wibisono]{duchi2015optimal}
Duchi, J.~C., Jordan, M.~I., Wainwright, M.~J., and Wibisono, A.
\newblock Optimal rates for zero-order convex optimization: The power of two
  function evaluations.
\newblock \emph{IEEE Transactions on Information Theory}, 61\penalty0
  (5):\penalty0 2788--2806, 2015.

\bibitem[Fang et~al.(2018)Fang, Li, Lin, and Zhang]{fang2018spider}
Fang, C., Li, C.~J., Lin, Z., and Zhang, T.
\newblock Spider: near-optimal non-convex optimization via stochastic path
  integrated differential estimator.
\newblock In \emph{Proceedings of the 32nd International Conference on Neural
  Information Processing Systems}, pp.\  687--697, 2018.

\bibitem[Gao \& Huang(2020)Gao and Huang]{gao2020can}
Gao, H. and Huang, H.
\newblock Can stochastic zeroth-order frank-wolfe method converge faster for
  non-convex problems?
\newblock In \emph{International Conference on Machine Learning}, pp.\
  3377--3386. PMLR, 2020.

\bibitem[Gao et~al.(2018)Gao, Jiang, and Zhang]{gao2018information}
Gao, X., Jiang, B., and Zhang, S.
\newblock On the information-adaptive variants of the admm: an iteration
  complexity perspective.
\newblock \emph{Journal of Scientific Computing}, 76\penalty0 (1):\penalty0
  327--363, 2018.

\bibitem[Ghadimi \& Lan(2013)Ghadimi and Lan]{ghadimi2013stochastic}
Ghadimi, S. and Lan, G.
\newblock Stochastic first-and zeroth-order methods for nonconvex stochastic
  programming.
\newblock \emph{SIAM Journal on Optimization}, 23\penalty0 (4):\penalty0
  2341--2368, 2013.

\bibitem[Ghadimi et~al.(2016)Ghadimi, Lan, and Zhang]{ghadimi2016mini}
Ghadimi, S., Lan, G., and Zhang, H.
\newblock Mini-batch stochastic approximation methods for nonconvex stochastic
  composite optimization.
\newblock \emph{Mathematical Programming}, 155\penalty0 (1-2):\penalty0
  267--305, 2016.

\bibitem[Gu et~al.(2016)Gu, Huo, and Huang]{gu2016zeroth}
Gu, B., Huo, Z., and Huang, H.
\newblock Zeroth-order asynchronous doubly stochastic algorithm with variance
  reduction.
\newblock \emph{arXiv preprint arXiv:1612.01425}, 2016.

\bibitem[Gu et~al.(2021{\natexlab{a}})Gu, Liu, Zhang, Geng, and
  Huang]{gu2021optimizing}
Gu, B., Liu, G., Zhang, Y., Geng, X., and Huang, H.
\newblock Optimizing large-scale hyperparameters via automated learning
  algorithm.
\newblock \emph{arXiv preprint arXiv:2102.09026}, 2021{\natexlab{a}}.

\bibitem[Gu et~al.(2021{\natexlab{b}})Gu, Wei, Gao, Xiong, Deng, and
  Huang]{gu2021black}
Gu, B., Wei, X., Gao, S., Xiong, Z., Deng, C., and Huang, H.
\newblock Black-box reductions for zeroth-order gradient algorithms to achieve
  lower query complexity.
\newblock \emph{Journal of Machine Learning Research}, 22\penalty0
  (170):\penalty0 1--47, 2021{\natexlab{b}}.

\bibitem[Guo et~al.(2021)Guo, Xu, Yin, Jin, and Yang]{guo2021stochastic}
Guo, Z., Xu, Y., Yin, W., Jin, R., and Yang, T.
\newblock On stochastic moving-average estimators for non-convex optimization.
\newblock \emph{arXiv preprint arXiv:2104.14840}, 2021.

\bibitem[Huang \& Huang(2021)Huang and Huang]{huang2021biadam}
Huang, F. and Huang, H.
\newblock Biadam: Fast adaptive bilevel optimization methods.
\newblock \emph{arXiv preprint arXiv:2106.11396}, 2021.

\bibitem[Huang et~al.(2020{\natexlab{a}})Huang, Gao, Pei, and
  Huang]{huang2020accelerated}
Huang, F., Gao, S., Pei, J., and Huang, H.
\newblock Accelerated zeroth-order and first-order momentum methods from mini
  to minimax optimization.
\newblock \emph{arXiv preprint arXiv:2008.08170}, 2020{\natexlab{a}}.

\bibitem[Huang et~al.(2020{\natexlab{b}})Huang, Tao, and
  Chen]{huang2020accelerated2}
Huang, F., Tao, L., and Chen, S.
\newblock Accelerated stochastic gradient-free and projection-free methods.
\newblock In \emph{International Conference on Machine Learning}, pp.\
  4519--4530. PMLR, 2020{\natexlab{b}}.

\bibitem[Huang et~al.(2021)Huang, Li, and Huang]{huang2021super}
Huang, F., Li, J., and Huang, H.
\newblock Super-adam: faster and universal framework of adaptive gradients.
\newblock \emph{Advances in Neural Information Processing Systems}, 34, 2021.

\bibitem[Ji et~al.(2019)Ji, Wang, Zhou, and Liang]{ji2019improved}
Ji, K., Wang, Z., Zhou, Y., and Liang, Y.
\newblock Improved zeroth-order variance reduced algorithms and analysis for
  nonconvex optimization.
\newblock In \emph{International conference on machine learning}, pp.\
  3100--3109. PMLR, 2019.

\bibitem[Johnson \& Zhang(2013)Johnson and Zhang]{johnson2013accelerating}
Johnson, R. and Zhang, T.
\newblock Accelerating stochastic gradient descent using predictive variance
  reduction.
\newblock \emph{Advances in neural information processing systems},
  26:\penalty0 315--323, 2013.

\bibitem[Lei et~al.(2017)Lei, Ju, Chen, and Jordan]{lei2017non}
Lei, L., Ju, C., Chen, J., and Jordan, M.~I.
\newblock Non-convex finite-sum optimization via scsg methods.
\newblock In \emph{Proceedings of the 31st International Conference on Neural
  Information Processing Systems}, pp.\  2345--2355, 2017.

\bibitem[Lin et~al.(2014)Lin, Lu, and Xiao]{lin2014accelerated}
Lin, Q., Lu, Z., and Xiao, L.
\newblock An accelerated proximal coordinate gradient method.
\newblock \emph{Advances in Neural Information Processing Systems},
  27:\penalty0 3059--3067, 2014.

\bibitem[Lin et~al.(2019)Lin, Ma, and Xu]{lin2019inexact}
Lin, Q., Ma, R., and Xu, Y.
\newblock Inexact proximal-point penalty methods for constrained non-convex
  optimization.
\newblock \emph{arXiv preprint arXiv:1908.11518}, 2019.

\bibitem[Lin et~al.(2020)Lin, Jin, and Jordan]{lin2020gradient}
Lin, T., Jin, C., and Jordan, M.
\newblock On gradient descent ascent for nonconvex-concave minimax problems.
\newblock In \emph{International Conference on Machine Learning}, pp.\
  6083--6093. PMLR, 2020.

\bibitem[Liu et~al.(2020)Liu, Zhang, Orabona, and Yang]{liu2020adam}
Liu, M., Zhang, W., Orabona, F., and Yang, T.
\newblock Adam $^+$: A stochastic method with adaptive variance reduction.
\newblock \emph{arXiv preprint arXiv:2011.11985}, 2020.

\bibitem[Liu et~al.(2018{\natexlab{a}})Liu, Chen, Chen, and
  Hero]{liu2018zeroth}
Liu, S., Chen, J., Chen, P.-Y., and Hero, A.
\newblock Zeroth-order online alternating direction method of multipliers:
  Convergence analysis and applications.
\newblock In \emph{International Conference on Artificial Intelligence and
  Statistics}, pp.\  288--297. PMLR, 2018{\natexlab{a}}.

\bibitem[Liu et~al.(2018{\natexlab{b}})Liu, Kailkhura, Chen, Ting, Chang, and
  Amini]{liu2018zeroth_var}
Liu, S., Kailkhura, B., Chen, P.-Y., Ting, P., Chang, S., and Amini, L.
\newblock Zeroth-order stochastic variance reduction for nonconvex
  optimization.
\newblock \emph{Advances in Neural Information Processing Systems},
  31:\penalty0 3727--3737, 2018{\natexlab{b}}.

\bibitem[Liu et~al.(2018{\natexlab{c}})Liu, Li, Chen, Haupt, and
  Amini]{liu2018zeroth_pro}
Liu, S., Li, X., Chen, P.-Y., Haupt, J., and Amini, L.
\newblock Zeroth-order stochastic projected gradient descent for nonconvex
  optimization.
\newblock In \emph{2018 IEEE Global Conference on Signal and Information
  Processing (GlobalSIP)}, pp.\  1179--1183. IEEE, 2018{\natexlab{c}}.

\bibitem[Nesterov \& Spokoiny(2017)Nesterov and Spokoiny]{nesterov2017random}
Nesterov, Y. and Spokoiny, V.
\newblock Random gradient-free minimization of convex functions.
\newblock \emph{Foundations of Computational Mathematics}, 17\penalty0
  (2):\penalty0 527--566, 2017.

\bibitem[Nguyen et~al.(2017)Nguyen, Liu, Scheinberg, and
  Tak{\'a}{\v{c}}]{nguyen2017sarah}
Nguyen, L.~M., Liu, J., Scheinberg, K., and Tak{\'a}{\v{c}}, M.
\newblock Sarah: A novel method for machine learning problems using stochastic
  recursive gradient.
\newblock In \emph{International Conference on Machine Learning}, pp.\
  2613--2621. PMLR, 2017.

\bibitem[Roux et~al.(2012)Roux, Schmidt, and Bach]{roux2012stochastic}
Roux, N.~L., Schmidt, M., and Bach, F.
\newblock A stochastic gradient method with an exponential convergence rate for
  finite training sets.
\newblock \emph{arXiv preprint arXiv:1202.6258}, 2012.

\bibitem[Shi \& Gu(2021)Shi and Gu]{shi2021improved}
Shi, W. and Gu, B.
\newblock Improved penalty method via doubly stochastic gradients for bilevel
  hyperparameter optimization.
\newblock In \emph{Proceedings of the AAAI Conference on Artificial
  Intelligence}, volume~35, pp.\  9621--9629, 2021.

\bibitem[Tran-Dinh et~al.(2021)Tran-Dinh, Pham, Phan, and
  Nguyen]{tran2021hybrid}
Tran-Dinh, Q., Pham, N.~H., Phan, D.~T., and Nguyen, L.~M.
\newblock A hybrid stochastic optimization framework for composite nonconvex
  optimization.
\newblock \emph{Mathematical Programming}, pp.\  1--67, 2021.

\bibitem[Wang et~al.(2017)Wang, Fang, and Liu]{wang2017stochastic}
Wang, M., Fang, E.~X., and Liu, H.
\newblock Stochastic compositional gradient descent: algorithms for minimizing
  compositions of expected-value functions.
\newblock \emph{Mathematical Programming}, 161\penalty0 (1-2):\penalty0
  419--449, 2017.

\bibitem[Wang et~al.(2018)Wang, Du, Balakrishnan, and
  Singh]{wang2018stochastic}
Wang, Y., Du, S., Balakrishnan, S., and Singh, A.
\newblock Stochastic zeroth-order optimization in high dimensions.
\newblock In \emph{International Conference on Artificial Intelligence and
  Statistics}, pp.\  1356--1365. PMLR, 2018.

\bibitem[Wang et~al.(2019)Wang, Ji, Zhou, Liang, and
  Tarokh]{wang2019spiderboost}
Wang, Z., Ji, K., Zhou, Y., Liang, Y., and Tarokh, V.
\newblock Spiderboost and momentum: Faster variance reduction algorithms.
\newblock \emph{Advances in Neural Information Processing Systems},
  32:\penalty0 2406--2416, 2019.

\bibitem[Wang et~al.(2023)Wang, Balasubramanian, Ma, and
  Razaviyayn]{wang2020zeroth}
Wang, Z., Balasubramanian, K., Ma, S., and Razaviyayn, M.
\newblock Zeroth-order algorithms for nonconvex--strongly-concave minimax
  problems with improved complexities.
\newblock \emph{Journal of Global Optimization}, 87\penalty0 (2):\penalty0
  709--740, 2023.

\bibitem[Wei et~al.(2021)Wei, Gu, and Huang]{wei2021accelerated}
Wei, X., Gu, B., and Huang, H.
\newblock An accelerated variance-reduced conditional gradient sliding
  algorithm for first-order and zeroth-order optimization.
\newblock \emph{arXiv preprint arXiv:2109.08858}, 2021.

\bibitem[Zafar et~al.(2017)Zafar, Valera, Gomez~Rodriguez, and
  Gummadi]{zafar2017fairness}
Zafar, M.~B., Valera, I., Gomez~Rodriguez, M., and Gummadi, K.~P.
\newblock Fairness beyond disparate treatment \& disparate impact: Learning
  classification without disparate mistreatment.
\newblock In \emph{Proceedings of the 26th international conference on world
  wide web}, pp.\  1171--1180, 2017.

\bibitem[Zhou et~al.(2018)Zhou, Xu, and Gu]{zhou2018stochastic}
Zhou, D., Xu, P., and Gu, Q.
\newblock Stochastic nested variance reduction for nonconvex optimization.
\newblock In \emph{Proceedings of the 32nd International Conference on Neural
  Information Processing Systems}, pp.\  3925--3936, 2018.

\end{thebibliography}
\bibliographystyle{icml2022}
\appendix
\newpage
\onecolumn
\section{Detailed Proofs}
\subsection{Proof of Proposition \ref{prop1}}
\begin{proof}
	Since the $(\boldsymbol{w}^*,\boldsymbol{p}^*)$ is the $\epsilon$-stationary point of $\min_{\boldsymbol{w}}\max_{\boldsymbol{p}\in\Delta^m}\mathcal{L}(\boldsymbol{w},\boldsymbol{p})$, then we have 
	\begin{align}
	\|\nabla_{\boldsymbol{w}}f_0(\boldsymbol{w}^*) + \beta\sum_{j=1}^mp_j^*2\max\{f_j(\boldsymbol{w}^*),0\}\nabla_{\boldsymbol{w}}f_j(\boldsymbol{w}^*)\|_2^2\leq\epsilon^2.
	\end{align}
	Let $\alpha_j^*=2\beta p_j^*\max\{f_j(\boldsymbol{w}^*),0\}$ and $\epsilon\leq\epsilon_1$, we have 
	\begin{align}
	\|\nabla_{\boldsymbol{w}}f_0(\boldsymbol{w}^*) + \sum_{j=1}^m\alpha_j^*\nabla_{\boldsymbol{w}}f_j(\boldsymbol{w}^*)\|_2^2\leq\epsilon_1^2.
	\end{align}
	Then the first condition in Definition 2 is satisfied. 
	
	Using $\|\nabla_{\boldsymbol{p}}\mathcal{L}(\boldsymbol{w}^*,\boldsymbol{p}^*)\|_2^2\leq\epsilon^2$ and $0\leq p_j^2\leq 1$, we have 
	\begin{align}
	\sum_{j=1}^m (\beta \phi_j(\boldsymbol{w}^*)-\lambda p_j^*)^2\leq \epsilon^2.
	\end{align}
	Using the inequality $\|a+b\|_2^2\leq2\|a\|_2^2+2\|b\|_2^2$, we have
	\begin{align}
	&\dfrac{1}{2}\beta^2\sum_{j=1}^{m}\phi_j(\boldsymbol{w}^*)^2\nonumber\\
	\leq&
	\sum_{j=1}^{m}(\beta\phi_j(\boldsymbol{w}^*)-\lambda p_j^*)^2+\lambda^2\sum_{j=1}^m(p_j^*)^2\nonumber\\
	\leq& \epsilon^2+m\lambda^2.
	\end{align}
	Then, using $(\dfrac{\sum_{i=0}^na_i}{n})^2\leq\dfrac{\sum_{i=0}^na_i^2}{n}$, we have 
	\begin{align}
	\sum_{j=1}^m\phi_j(\boldsymbol{w}^*)\leq \sqrt{m\sum_{j=1}^m\phi_j(\boldsymbol{w}^*)^2}\leq\sqrt{\dfrac{2m\epsilon^2+2m^2\lambda^2}{\beta^2}}.
	\end{align}
	Let $\sqrt{\dfrac{2m\epsilon^2+2m^2\lambda^2}{\beta^2}}\leq\epsilon_2^2$ and $\phi_j(\boldsymbol{w})=(\max\{f_j(\boldsymbol{w}),0\})^2$, we can obtain
	\begin{align}\label{minmax_to_ekkt_2}
	\sum_{j=1}^m (\max\{f_j(\boldsymbol{w}^*),0\})^2\leq\epsilon_2^2.
	\end{align}
	Therefore, the second condition in Definition 2 is satisfied.
	
	Based on the inequality $\|\langle\boldsymbol{a},\boldsymbol{b}\rangle\|_2^2\leq\|\boldsymbol{a}\|_2^2\|\boldsymbol{b}\|_2^2$, we can multiply $\sum_{j=1}^m(\alpha_j^*)^2$ on both sides of the inequality \ref{minmax_to_ekkt_2}, such that we have	
	\begin{align}
	\left(\sum_{j=1}^m \alpha_j^*\max\{f_j(\boldsymbol{w}^*),0\}\right)^2\leq\sum_{j=1}^m (\alpha_j^*)^2	\sum_{j=1}^m \max\{f_j(\boldsymbol{w}^*),0\}^2\leq\epsilon_2^2\sum_{j=1}^m(\alpha_j^*)^2.
	\end{align}
	Since $\alpha_j^*\geq0$ and $\max\{f_j(\boldsymbol{w}^*),0\}\geq 0$
	\begin{align}
	\sum_{j=1}^m (\alpha_j^*\max\{f_j(\boldsymbol{w}^*),0\})^2\leq\left(\sum_{j=1}^m \alpha_j^*\max\{f_j(\boldsymbol{w}^*),0\}\right)^2\leq\epsilon_2^2\sum_{j=1}^m(\alpha_j^*)^2.
	\end{align}
	Using inequality \ref{minmax_to_ekkt_2}, we have  $(\alpha_j^*)^2=4\beta^2 (p_j^*)^2(\max\{f_j(\boldsymbol{w}^*),0\})^2\leq4\beta^2\epsilon_2^2$,
	Let $4\beta^2\epsilon_2^2\leq\epsilon_3^2$, we have
	\begin{align}
	\sum_{j=1}^m (\alpha_j^*\max\{f_j(\boldsymbol{w}^*),0\})^2\leq\epsilon_3^2.
	\end{align} 
	If $f_j(\boldsymbol{w}^*)\leq0$, we have $\alpha_j^*=2\beta p_j^*\max\{f_j(\boldsymbol{w}^*),0\}=0$. Therefore, we have 
	\begin{align}
	\sum_{j=1}^m (\alpha_j^*f_j(\boldsymbol{w}^*))^2\leq\epsilon_3^2,
	\end{align}
	which means that the third condition in Definition 2 is satisfied.
	
	That completes the proof.
\end{proof}

\subsection{Proof of Proposition \ref{prop2}}

\begin{proof}
	Assume that a point $\hat{\boldsymbol{w}}$ satisfies that $\|\nabla_{\boldsymbol{w}} g(\hat{\boldsymbol{w}})\|_2\leq \epsilon$, the optimization problem $\max_{\boldsymbol{p}\in\Delta^m}\mathcal{L}(\hat{\boldsymbol{w}},\boldsymbol{p})$ is strongly concave w.r.t $\boldsymbol{p}$ and $\boldsymbol{p}^*(\hat{\boldsymbol{w}})$ is uniquely defined. Solving this this strongly concave problem $\max_{\boldsymbol{p}\in\Delta^m}\mathcal{L}(\hat{\boldsymbol{w}},\boldsymbol{p})$, we can obtain a point $\boldsymbol{p}'$ satisfying that 
	\begin{align}
	\|\nabla_{\boldsymbol{p}}\mathcal{L}(\hat{\boldsymbol{w}},\boldsymbol{p}')\|_2\leq\epsilon \ and \ \|\boldsymbol{p}'-\boldsymbol{p}^*(\hat{\boldsymbol{w}})\|_2\leq\epsilon.
	\end{align}
	
	If $\|\nabla_{\boldsymbol{w}}g(\hat{\boldsymbol{w}})\|_2\leq\epsilon$, we have 
	\begin{align}
	&\|\nabla_{\boldsymbol{w}}\mathcal{L}(\hat{\boldsymbol{w}},\boldsymbol{p}')\|_2\nonumber\\
	\leq&\|\nabla_{\boldsymbol{w}}\mathcal{L}(\hat{\boldsymbol{w}},\boldsymbol{p}')-\nabla_{\boldsymbol{w}}g(\hat{\boldsymbol{w}})\|_2+\|\nabla_{\boldsymbol{w}}g(\hat{\boldsymbol{w}})\|_2\nonumber\\
	=&\|\nabla_{\boldsymbol{w}}\mathcal{L}(\hat{\boldsymbol{w}},\boldsymbol{p}')-\nabla_{\boldsymbol{w}}\mathcal{L}(\hat{\boldsymbol{w}},\boldsymbol{p}^*(\hat{\boldsymbol{w}}))\|_2+\epsilon\nonumber\\
	\leq&L\|\boldsymbol{p}'-\boldsymbol{p}^*(\hat{\boldsymbol{w}})\|_2+\epsilon\nonumber\\
	=&\mathcal{O}(\epsilon).
	\end{align}
\end{proof}

\subsection{Proof of Lemma \ref{lemma:bound_of_g_in_ASZOG}}

\begin{proof}
    Under Assumptions \ref{assum:penalty_function} and \ref{assum:bound_step}, we have  
	\begin{align}
	&g(\boldsymbol{w}_{t+1})\nonumber\\
	\leq&g(\boldsymbol{w}_t)+\nabla g(\boldsymbol{w}_t)^T(\boldsymbol{w}_{t+1}-\boldsymbol{w}_t)+\dfrac{L}{2}\|\boldsymbol{w}_{t+1}-\boldsymbol{w}_t\|_2^2\nonumber\\
	=&g(\boldsymbol{w}_t)-\eta_w\nabla g(\boldsymbol{w}_t)^T\dfrac{\boldsymbol{z}_{\boldsymbol{w}}^t}{\sqrt{\|\boldsymbol{z}_{\boldsymbol{w}}^t\|_2}+c}+\dfrac{L}{2}\dfrac{\|\boldsymbol{z}_{\boldsymbol{w}}^t\|_2^2}{\|\sqrt{\|\boldsymbol{z}_{\boldsymbol{w}}^t\|_2}+c\|_2^2}\nonumber\\
	\leq&g(\boldsymbol{w}_t)-\eta_wc_{1,l}\nabla g(\boldsymbol{w}_t)^T\boldsymbol{z}_{\boldsymbol{w}}^t+\dfrac{L\eta_w^2c_{1,u}^2}{2}\|\boldsymbol{z}_{\boldsymbol{w}}^t\|_2^2\nonumber\\
	=&g(\boldsymbol{w}_t)+\dfrac{\eta_wc_{1,l}}{2}\|\nabla g(\boldsymbol{w}_t)-\boldsymbol{z}_{\boldsymbol{w}}^t\|_2^2-\dfrac{\eta_wc_{1,l}}{2}\|\nabla g(\boldsymbol{w}_t)\|_2^2-\dfrac{\eta_wc_{1,l}}{2}\|\boldsymbol{z}_{\boldsymbol{w}}^t\|_2^2+\dfrac{L\eta_w^2c_{1,u}^2}{2}\|\boldsymbol{z}_{\boldsymbol{w}}^t\|_2^2\nonumber\\
	=&g(\boldsymbol{w}_t)-\dfrac{\eta_wc_{1,l}}{2}\|\nabla g(\boldsymbol{w}_t)\|_2^2-\dfrac{\eta_wc_{1,l}}{2}\|\boldsymbol{z}_{\boldsymbol{w}}^t\|_2^2+\dfrac{L\eta_w^2c_{1,u}^2}{2}\|\boldsymbol{z}_{\boldsymbol{w}}^t\|_2^2\nonumber\\
	&+\dfrac{\eta_wc_{1,l}}{2}\|\nabla g(\boldsymbol{w}_t)-\nabla g_{\mu}(\boldsymbol{w}_t)+\nabla g_{\mu}(\boldsymbol{w}_t)-\nabla_{\boldsymbol{w}}\mathcal{L}_{\mu}(\boldsymbol{w}_t,\boldsymbol{p}_t)+\nabla_{\boldsymbol{w}}\mathcal{L}_{\mu}(\boldsymbol{w}_t,\boldsymbol{p}_t)-\nabla_{\boldsymbol{w}}\mathcal{L}(\boldsymbol{w}_t,\boldsymbol{p}_t)\nonumber\\
	&+\nabla_{\boldsymbol{w}}\mathcal{L}(\boldsymbol{w}_t,\boldsymbol{p}_t)-\boldsymbol{z}_{\boldsymbol{w}}^t\|_2^2\nonumber\\
	\leq&g(\boldsymbol{w}_t)-\dfrac{\eta_wc_{1,l}}{2}\|\nabla g(\boldsymbol{w}_t)\|_2^2-\dfrac{\eta_wc_{1,l}}{2}\|\boldsymbol{z}_{\boldsymbol{w}}^t\|_2^2+\dfrac{L\eta_w^2c_{1,u}^2}{2}\|\boldsymbol{z}_{\boldsymbol{w}}^t\|_2^2\nonumber\\
	&+\eta_wc_{1,l}\|\nabla_{\boldsymbol{w}}\mathcal{L}(\boldsymbol{w}_t,\boldsymbol{p}^*(\boldsymbol{w}_t))-\nabla_{\boldsymbol{w}}\mathcal{L}_{\mu}(\boldsymbol{w}_t,\boldsymbol{p}^*(\boldsymbol{w}_t))\|_2^2+\eta_wc_{1,l}\|\nabla \nabla_{\boldsymbol{w}}\mathcal{L}_{\mu}(\boldsymbol{w}_t,\boldsymbol{p}^*(\boldsymbol{w}_t))-\nabla_{\boldsymbol{w}}\mathcal{L}_{\mu}(\boldsymbol{w}_t,\boldsymbol{p}_t)\|_2^2\nonumber\\
	&+\eta_wc_{1,l}\|\nabla_{\boldsymbol{w}}\mathcal{L}_{\mu}(\boldsymbol{w}_t,\boldsymbol{p}_t)-\nabla_{\boldsymbol{w}}\mathcal{L}(\boldsymbol{w}_t,\boldsymbol{p}_t)\|_2^2+\eta_wc_{1,l}\|\nabla_{\boldsymbol{w}}\mathcal{L}(\boldsymbol{w}_t,\boldsymbol{p}_t)-\boldsymbol{z}_{\boldsymbol{w}}^t\|_2^2\nonumber\\
	\leq&g(\boldsymbol{w}_t)-\dfrac{\eta_wc_{1,l}}{2}\|\nabla g(\boldsymbol{w}_t)\|_2^2-\dfrac{\eta_wc_{1,l}}{2}\|\boldsymbol{z}_{\boldsymbol{w}}^t\|_2^2+\dfrac{L\eta_w^2c_{1,u}^2}{2}\|\boldsymbol{z}_{\boldsymbol{w}}^t\|_2^2+\eta_wc_{1,l}\dfrac{\mu^2 L^2(d+3)^3}{4}+\eta_wc_{1,l}L^2\|\boldsymbol{p}^*(\boldsymbol{w}_t)-\boldsymbol{p}_t\|_2^2\nonumber\\
	&+\eta_wc_{1,l}\dfrac{\mu^2 L^2(d+3)^3}{4}+\eta_wc_{1,l}\|\nabla_{\boldsymbol{w}}\mathcal{L}(\boldsymbol{w}_t,\boldsymbol{p}_t)-\boldsymbol{z}_{\boldsymbol{w}}^t\|_2^2\nonumber\\
	\leq&g(\boldsymbol{w}_t)-\dfrac{\eta_wc_{1,l}}{2}\|\nabla g(\boldsymbol{w}_t)\|_2^2-\dfrac{\eta_wc_{1,l}}{2}\|\boldsymbol{z}_{\boldsymbol{w}}^t\|_2^2+\dfrac{L\eta_w^2c_{1,u}^2}{2}\|\boldsymbol{z}_{\boldsymbol{w}}^t\|_2^2+\eta_wc_{1,l}\dfrac{\mu^2 L^2(d+3)^3}{2}+\eta_wc_{1,l}L^2\|\boldsymbol{p}^*(\boldsymbol{w}_t)-\boldsymbol{p}_t\|_2^2\nonumber\\
	&+\eta_wc_{1,l}\|\nabla_{\boldsymbol{w}}\mathcal{L}(\boldsymbol{w}_t,\boldsymbol{p}_t)-\boldsymbol{z}_{\boldsymbol{w}}^t\|_2^2\nonumber\\
	\leq&g(\boldsymbol{w}_t)-\dfrac{\eta_wc_{1,l}}{2}\|\nabla g(\boldsymbol{w}_t)\|_2^2-\dfrac{\eta_wc_{1,l}}{4}\|\boldsymbol{z}_{\boldsymbol{w}}^t\|_2^2+\eta_wc_{1,l}\dfrac{\mu^2 L^2(d+3)^3}{2}+\eta_wc_{1,l}L^2\|\boldsymbol{p}^*(\boldsymbol{w}_t)-\boldsymbol{p}_t\|_2^2\nonumber\\
	&+\eta_wc_{1,l}\|\nabla_{\boldsymbol{w}}\mathcal{L}(\boldsymbol{w}_t,\boldsymbol{p}_t)-\boldsymbol{z}_{\boldsymbol{w}}^t\|_2^2.\nonumber\\
	\end{align}
	The last inequality is due to $\eta_{w}L\leq\dfrac{c_{1,l}}{2c_{1,u}^2}$.
\end{proof}
\subsection{Proof of Lemma \ref{lemma:bound_of_p_in_ASZOGD}}
\begin{proof}
	According to the update rule of $\boldsymbol{p}$, we have 
	\begin{align}
    	&\|\boldsymbol{p}_{t+1}-\boldsymbol{p}^*(\boldsymbol{w}_t)\|_2^2\nonumber\\
    	\leq&\|(1-a)\boldsymbol{p}_{t}+a\hat{\boldsymbol{p}}_{t+1}-\boldsymbol{p}^*(\boldsymbol{w}_t)\|_2^2\nonumber\\
    	=&\|\boldsymbol{p}_t-\boldsymbol{p}^*(\boldsymbol{w}_t)\|_2^2+a^2\|\boldsymbol{p}_t-\hat{\boldsymbol{p}}_{t+1}\|_2^2+2a\langle\boldsymbol{p}_t-\boldsymbol{p}^*(\boldsymbol{w}_t), \boldsymbol{p}_t-\hat{\boldsymbol{p}}_{t+1} \rangle.
	\end{align}
	Rearrange the above inequality, we have 
	\begin{align}
	\langle\boldsymbol{p}_t-\boldsymbol{p}^*(\boldsymbol{w}_t), \boldsymbol{p}_t-\hat{\boldsymbol{p}}_{t+1}  \rangle \geq \dfrac{1}{2a}\left( \|\boldsymbol{p}_{t+1}-\boldsymbol{p}^*(\boldsymbol{w}_t)\|_2^2-\|\boldsymbol{p}_t-\boldsymbol{p}^*(\boldsymbol{w}_t)\|_2^2 -a^2\|\boldsymbol{p}_t-\hat{\boldsymbol{p}}_{t+1}\|_2^2\right).
	\end{align}
	
	Due to the Assumption \ref{assum:penalty_function}, we have 
	\begin{align}
	\mathcal{L}(\boldsymbol{w}_t,\hat{\boldsymbol{p}}_{t+1})\geq\mathcal{L}(\boldsymbol{w}_t,\boldsymbol{p}_t)+\nabla_{\boldsymbol{p}}\mathcal{L}(\boldsymbol{w}_t,\boldsymbol{p}_t)^T(\hat{\boldsymbol{p}}_{t+1}-\boldsymbol{p}_t)-\dfrac{L}{2}\|\hat{\boldsymbol{p}}_{t+1}-\boldsymbol{p}_t\|_2^2.
	\end{align}
	
	In addition, according to the strongly concave, we have 
	\begin{align}
	&\mathcal{L}(\boldsymbol{w}_t,\boldsymbol{p})\nonumber\\
	\leq&
	\mathcal{L}(\boldsymbol{w}_t,\boldsymbol{p}_t)+\nabla_{\boldsymbol{p}}\mathcal{L}(\boldsymbol{w}_t,\boldsymbol{p}_t)^T(\boldsymbol{p}-\boldsymbol{p}_t)-\dfrac{\tau}{2}\|\boldsymbol{p}-\boldsymbol{p}_t\|_2^2\nonumber\\
	=&\mathcal{L}(\boldsymbol{w}_t,\boldsymbol{p}_t)+\nabla_{\boldsymbol{p}}\mathcal{L}(\boldsymbol{w}_t,\boldsymbol{p}_t)^T(\boldsymbol{p}-\hat{\boldsymbol{p}}_{t+1}+\hat{\boldsymbol{p}}_{t+1}-\boldsymbol{p}_t)-\dfrac{\tau}{2}\|\boldsymbol{p}-\boldsymbol{p}_t\|_2^2\nonumber\\
	=&\mathcal{L}(\boldsymbol{w}_t,\boldsymbol{p}_t)+\nabla_{\boldsymbol{p}}\mathcal{L}(\boldsymbol{w}_t,\boldsymbol{p}_t)^T(\boldsymbol{p}-\hat{\boldsymbol{p}}_{t+1})+\nabla_{\boldsymbol{p}}\mathcal{L}(\boldsymbol{w}_t,\boldsymbol{p}_t)^T(\hat{\boldsymbol{p}}_{t+1}-\boldsymbol{p}_t)-\dfrac{\tau}{2}\|\boldsymbol{p}-\boldsymbol{p}_t\|_2^2\nonumber\\
	=&\mathcal{L}(\boldsymbol{w}_t,\boldsymbol{p}_t)+(\nabla_{\boldsymbol{p}}\mathcal{L}(\boldsymbol{w}_t,\boldsymbol{p}_t)-\boldsymbol{z}_{\boldsymbol{p}}^t)^T(\boldsymbol{p}-\hat{\boldsymbol{p}}_{t+1})+\langle\boldsymbol{z}_{\boldsymbol{p}}^t,\boldsymbol{p}-\hat{\boldsymbol{p}}_{t+1}\rangle+\nabla_{\boldsymbol{p}}\mathcal{L}(\boldsymbol{w}_t,\boldsymbol{p}_t)^T(\hat{\boldsymbol{p}}_{t+1}-\boldsymbol{p}_t)-\dfrac{\tau}{2}\|\boldsymbol{p}-\boldsymbol{p}_t\|_2^2.
	\end{align}
	Then, using the above inequalities, we have 
	\begin{align}
	&\mathcal{L}(\boldsymbol{w}_t,\boldsymbol{p})\nonumber\\
	\leq&\mathcal{L}(\boldsymbol{w}_t,\hat{\boldsymbol{p}}_{t+1})+(\nabla_{\boldsymbol{p}}\mathcal{L}(\boldsymbol{w}_t,\boldsymbol{p}_t)-\boldsymbol{z}_{\boldsymbol{p}}^t)^T(\boldsymbol{p}-\hat{\boldsymbol{p}}_{t+1})+\langle\boldsymbol{z}_{\boldsymbol{p}}^t,\boldsymbol{p}-\hat{\boldsymbol{p}}_{t+1}\rangle-\dfrac{\tau}{2}\|\boldsymbol{p}-\boldsymbol{p}_t\|_2^2+\dfrac{L}{2}\|\hat{\boldsymbol{p}}_{t+1}-\boldsymbol{p}_t\|_2^2.
	\end{align}
	
	Due to the update rule of $\hat{\boldsymbol{p}}$, we have 
	\begin{align}
	\langle \hat{\boldsymbol{p}}_{t+1}-\boldsymbol{p}_t-\eta_{\boldsymbol{p}}\dfrac{\boldsymbol{z}_{\boldsymbol{p}}^t}{\sqrt{\|\boldsymbol{z}_{\boldsymbol{p}}^t\|_2}+c},\boldsymbol{p}-\hat{\boldsymbol{p}}_{t+1} \rangle\geq 0, \quad\forall\boldsymbol{p}\in\Delta^m.
	\end{align}
	Then, we have 
	\begin{align}
    	&\eta_{\boldsymbol{p}}c_{2,l}\langle \boldsymbol{z}_{\boldsymbol{p}}^t,\boldsymbol{p}-\hat{\boldsymbol{p}}_{t+1} \rangle\nonumber\\
    	\leq&\langle \eta_{\boldsymbol{p}}\dfrac{\boldsymbol{z}_{\boldsymbol{p}}^t}{\sqrt{\|\boldsymbol{z}_{\boldsymbol{p}}^t\|_2}+c},\boldsymbol{p}-\hat{\boldsymbol{p}}_{t+1} \rangle\nonumber\\
    	\leq&\langle \hat{\boldsymbol{p}}_{t+1}-\boldsymbol{p}_t,\boldsymbol{p}-\hat{\boldsymbol{p}}_{t+1} \rangle\nonumber\\
    	=&\langle \hat{\boldsymbol{p}}_{t+1}-\boldsymbol{p}_t,\boldsymbol{p}-\boldsymbol{p}_t+\boldsymbol{p}_t-\hat{\boldsymbol{p}}_{t+1} \rangle\nonumber\\
    	=&-\| \hat{\boldsymbol{p}}_{t+1}-\boldsymbol{p}_t \|_2^2+\langle \hat{\boldsymbol{p}}_{t+1}-\boldsymbol{p}_t,\boldsymbol{p}-\boldsymbol{p}_t\rangle.
	\end{align}
	
	In addition, we have 
	\begin{align}
	&\langle\nabla_{\boldsymbol{p}}\mathcal{L}(\boldsymbol{w}_t,\boldsymbol{p}_t)-\boldsymbol{z}_{\boldsymbol{p}}^t,\boldsymbol{p}^*(\boldsymbol{w}_{t})-\hat{\boldsymbol{p}}_{t+1}\rangle\nonumber\\
	=&\langle\nabla_{\boldsymbol{p}}\mathcal{L}(\boldsymbol{w}_t,\boldsymbol{p}_t)-\boldsymbol{z}_{\boldsymbol{p}}^t,\boldsymbol{p}^*(\boldsymbol{w}_{t})-\boldsymbol{p}_t\rangle + \langle\nabla_{\boldsymbol{p}}\mathcal{L}(\boldsymbol{w}_t,\boldsymbol{p}_t)-\boldsymbol{z}_{\boldsymbol{p}}^t,\boldsymbol{p}_t-\hat{\boldsymbol{p}}_{t+1}\rangle\nonumber\\
	\leq&\dfrac{1}{\tau}\|\nabla_{\boldsymbol{p}}\mathcal{L}(\boldsymbol{w}_t,\boldsymbol{p}_t)-\boldsymbol{z}_{\boldsymbol{p}}^t\|_2^2+\dfrac{\tau}{4}\|\boldsymbol{p}^*(\boldsymbol{w}_{t})-\boldsymbol{p}_t\|_2^2 +\dfrac{1}{\tau}\| \nabla_{\boldsymbol{p}}\mathcal{L}(\boldsymbol{w}_t,\boldsymbol{p}_t)-\boldsymbol{z}_{\boldsymbol{p}}^t\|_2^2+\dfrac{\tau}{4}\|\boldsymbol{p}_t-\hat{\boldsymbol{p}}_{t+1}\|_2^2\nonumber\\
	\leq&\dfrac{2}{\tau}\|\nabla_{\boldsymbol{p}}\mathcal{L}(\boldsymbol{w}_t,\boldsymbol{p}_t)-\boldsymbol{z}_{\boldsymbol{p}}^t\|_2^2+\dfrac{\tau}{4}\|\boldsymbol{p}^*(\boldsymbol{w}_{t})-\boldsymbol{p}_t\|_2^2 +\dfrac{\tau}{4}\|\boldsymbol{p}_t-\hat{\boldsymbol{p}}_{t+1}\|_2^2.
	\end{align}
	
	Then, we have 
	\begin{align}
	&\mathcal{L}(\boldsymbol{w}_t,\boldsymbol{p})\nonumber\\
	\leq&\mathcal{L}(\boldsymbol{w}_t,\hat{\boldsymbol{p}}_{t+1})-\dfrac{1}{\eta_{p}c_{2,l}}\| \hat{\boldsymbol{p}}_{t+1}-\boldsymbol{p}_t \|_2^2+\dfrac{1}{\eta_{p}c_{2,l}}\langle \hat{\boldsymbol{p}}_{t+1}-\boldsymbol{p}_t,\boldsymbol{p}-\boldsymbol{p}_t\rangle-\dfrac{\tau}{2}\|\boldsymbol{p}-\boldsymbol{p}_t\|_2^2+\dfrac{L}{2}\|\hat{\boldsymbol{p}}_{t+1}-\boldsymbol{p}_t\|_2^2\nonumber\\
	&+\dfrac{2}{\tau}\|\nabla_{\boldsymbol{p}}\mathcal{L}(\boldsymbol{w}_t,\boldsymbol{p}_t)-\boldsymbol{z}_{\boldsymbol{p}}^t\|_2^2+\dfrac{\tau}{4}\|\boldsymbol{p}^*(\boldsymbol{w}_{t})-\boldsymbol{p}_t\|_2^2 +\dfrac{\tau}{4}\|\boldsymbol{p}_t-\hat{\boldsymbol{p}}_{t+1}\|_2^2.
	\end{align}
	Let $\boldsymbol{p}=\boldsymbol{p}^*(\boldsymbol{w}_{t})$, we have
	\begin{align}
	&\mathcal{L}(\boldsymbol{w}_t,\hat{\boldsymbol{p}}_{t+1})\nonumber\\
	\leq&\mathcal{L}(\boldsymbol{w}_t,\boldsymbol{p}^*(\boldsymbol{w}_{t}))\nonumber\\
	\leq&\mathcal{L}(\boldsymbol{w}_t,\hat{\boldsymbol{p}}_{t+1})-\dfrac{1}{\eta_{p}c_{2,l}}\| \hat{\boldsymbol{p}}_{t+1}-\boldsymbol{p}_t \|_2^2-\dfrac{\tau}{2}\|\boldsymbol{p}^*(\boldsymbol{w}_{t})-\boldsymbol{p}_t\|_2^2+\dfrac{L}{2}\|\hat{\boldsymbol{p}}_{t+1}-\boldsymbol{p}_t\|_2^2\nonumber\\
	&+\dfrac{2}{\tau}\|\nabla_{\boldsymbol{p}}\mathcal{L}(\boldsymbol{w}_t,\boldsymbol{p}_t)-\boldsymbol{z}_{\boldsymbol{p}}^t\|_2^2+\dfrac{\tau}{4}\|\boldsymbol{p}^*(\boldsymbol{w}_{t})-\boldsymbol{p}_t\|_2^2 +\dfrac{\tau}{4}\|\boldsymbol{p}_t-\hat{\boldsymbol{p}}_{t+1}\|_2^2\nonumber\\
	&-\dfrac{1}{2a\eta_{p}c_{2,l}}\left( \|\boldsymbol{p}_{t+1}-\boldsymbol{p}^*(\boldsymbol{w}_t)\|_2^2-\|\boldsymbol{p}_t-\boldsymbol{p}^*(\boldsymbol{w}_t)\|_2^2 -a^2\|\boldsymbol{p}_t-\hat{\boldsymbol{p}}_{t+1}\|_2^2\right).
	\end{align} 
	
	Rearrange the inequality, we have 
	\begin{align}
	&\|\boldsymbol{p}_{t+1}-\boldsymbol{p}^*(\boldsymbol{w}_t)\|_2^2\nonumber\\
	\leq&-2a\eta_{p}c_{2,l}\left(\dfrac{1}{\eta_{p}c_{2,l}}-\dfrac{L}{2}-\dfrac{\tau}{4}-\dfrac{1}{2b\eta_{p}c_{2,l}}\right)\|\boldsymbol{p}_t-\hat{\boldsymbol{p}}_{t+1}\|_2^2+\dfrac{4a\eta_pc_{2,l}}{\tau}\|\nabla_{\boldsymbol{p}}\mathcal{L}(\boldsymbol{w}_t,\boldsymbol{p}_t)-\boldsymbol{z}_{\boldsymbol{p}}^t\|_2^2\nonumber\\
	&+(1-\dfrac{\tau a \eta_{p}c_{2,l}}{2})\|\boldsymbol{p}^*(\boldsymbol{w}_{t})-\boldsymbol{p}_t\|_2^2\nonumber\\
	\leq&-2a\eta_{p}c_{2,l}\left(\dfrac{1}{2\eta_{p}c_{2,l}}-\dfrac{3L}{4}\right)\|\boldsymbol{p}_t-\hat{\boldsymbol{p}}_{t+1}\|_2^2+\dfrac{4a\eta_pc_{2,l}}{\tau}\|\nabla_{\boldsymbol{p}}\mathcal{L}(\boldsymbol{w}_t,\boldsymbol{p}_t)-\boldsymbol{z}_{\boldsymbol{p}}^t\|_2^2\nonumber\\
	&+(1-\dfrac{\tau a \eta_{p}c_{2,l}}{2})\|\boldsymbol{p}^*(\boldsymbol{w}_{t})-\boldsymbol{p}_t\|_2^2.
	\end{align}
	where we use $a\leq1$, $\tau\leq L$ and $\eta_{p}\leq\dfrac{1}{3c_{2.l}L}$.
	
	Then, we have 
	\begin{align}
	&\|\boldsymbol{p}_{t+1}-\boldsymbol{p}^*(\boldsymbol{w}_{t+1})\|_2^2\nonumber\\
	=&\|\boldsymbol{p}_{t+1}-\boldsymbol{p}^*(\boldsymbol{w}_{t})+\boldsymbol{p}^*(\boldsymbol{w}_{t})-\boldsymbol{p}^*(\boldsymbol{w}_{t+1})\|_2^2\nonumber\\
	\leq&(1+\dfrac{\tau a \eta_{p}c_{2,l}}{4})\|\boldsymbol{p}_{t+1}-\boldsymbol{p}^*(\boldsymbol{w}_{t})\|_2^2+(1+\dfrac{4}{\tau a \eta_{p}c_{2,l}})\|\boldsymbol{p}^*(\boldsymbol{w}_{t})-\boldsymbol{p}^*(\boldsymbol{w}_{t+1})\|_2^2\nonumber\\
	\leq&-2a\eta_{p}c_{2,l}(1+\dfrac{\tau a \eta_{p}c_{2,l}}{4})\left(\dfrac{1}{2\eta_{p}c_{2,l}}-\dfrac{3L}{4}\right)\|\boldsymbol{p}_t-\hat{\boldsymbol{p}}_{t+1}\|_2^2+\dfrac{4a\eta_pc_{2,l}}{\tau}(1+\dfrac{\tau a \eta_{p}c_{2,l}}{4})\|\nabla_{\boldsymbol{p}}\mathcal{L}(\boldsymbol{w}_t,\boldsymbol{p}_t)-\boldsymbol{z}_{\boldsymbol{p}}^t\|_2^2\nonumber\\
	&+(1-\dfrac{\tau a \eta_{p}c_{2,l}}{2})(1+\dfrac{\tau a \eta_{p}c_{2,l}}{4})\|\boldsymbol{p}^*(\boldsymbol{w}_{t})-\boldsymbol{p}_t\|_2^2+(1+\dfrac{4}{\tau a \eta_{p}c_{2,l}})L_g^2\|\boldsymbol{w}_{t}-\boldsymbol{w}_{t+1}\|_2^2\nonumber\\
	\leq&-2a\eta_{p}c_{2,l}(1+\dfrac{\tau a \eta_{p}c_{2,l}}{4})\left(\dfrac{1}{2\eta_{p}c_{2,l}}-\dfrac{3L}{4}\right)\|\boldsymbol{p}_t-\hat{\boldsymbol{p}}_{t+1}\|_2^2+\dfrac{8a\eta_pc_{2,l}}{\tau}\|\nabla_{\boldsymbol{p}}\mathcal{L}(\boldsymbol{w}_t,\boldsymbol{p}_t)-\boldsymbol{z}_{\boldsymbol{p}}^t\|_2^2\nonumber\\
	&+(1-\dfrac{\tau a \eta_{p}c_{2,l}}{4})\|\boldsymbol{p}^*(\boldsymbol{w}_{t})-\boldsymbol{p}_t\|_2^2+\dfrac{8L^2_g}{\tau a \eta_{p}c_{2,l}}\|\boldsymbol{w}_{t}-\boldsymbol{w}_{t+1}\|_2^2\nonumber\\
	\leq&-\dfrac{2\eta_{p}c_{2,l}}{a}(1+\dfrac{\tau a \eta_{p}c_{2,l}}{4})\left(\dfrac{1}{2\eta_{p}c_{2,l}}-\dfrac{3L}{4}\right)\|\boldsymbol{p}_t-{\boldsymbol{p}}_{t+1}\|_2^2+\dfrac{8a\eta_pc_{2,l}}{\tau}\|\nabla_{\boldsymbol{p}}\mathcal{L}(\boldsymbol{w}_t,\boldsymbol{p}_t)-\boldsymbol{z}_{\boldsymbol{p}}^t\|_2^2\nonumber\\
	&+(1-\dfrac{\tau a \eta_{p}c_{2,l}}{4})\|\boldsymbol{p}^*(\boldsymbol{w}_{t})-\boldsymbol{p}_t\|_2^2+\dfrac{8L^2_g}{\tau a \eta_{p}c_{2,l}}\|\boldsymbol{w}_{t}-\boldsymbol{w}_{t+1}\|_2^2\nonumber\\
	\leq&-\dfrac{1}{4a}\|\boldsymbol{p}_t-{\boldsymbol{p}}_{t+1}\|_2^2+\dfrac{8a\eta_pc_{2,l}}{\tau}\|\nabla_{\boldsymbol{p}}\mathcal{L}(\boldsymbol{w}_t,\boldsymbol{p}_t)-\boldsymbol{z}_{\boldsymbol{p}}^t\|_2^2\nonumber\\
	&+(1-\dfrac{\tau a \eta_{p}c_{2,l}}{4})\|\boldsymbol{p}^*(\boldsymbol{w}_{t})-\boldsymbol{p}_t\|_2^2+\dfrac{8L^2_g}{\tau a \eta_{p}c_{2,l}}\|\boldsymbol{w}_{t}-\boldsymbol{w}_{t+1}\|_2^2.
	\end{align}
\end{proof}

\subsection{Proof of Lemma \ref{lemma:bound_of_z_in_ASZOGD}}
\begin{proof}
According to the update rule of $z_{\boldsymbol{p}}$, we have
	\begin{align}
	&\boldsymbol{z}_{\boldsymbol{p}}^{t+1}-\boldsymbol{z}_{\boldsymbol{p}}^{t}
	=-b\boldsymbol{z}_{\boldsymbol{p}}^t+bH^{t+1}.
	\end{align}
	Then, we have 
	\begin{align}
	&\mathbb{E}[\| \nabla_{\boldsymbol{p}}\mathcal{L}(\boldsymbol{w}_{t+1},\boldsymbol{p}_{t+1}) -\boldsymbol{z}_{\boldsymbol{p}}^{t+1} \|_2^2]\nonumber\\
	=&\mathbb{E}[\| \nabla_{\boldsymbol{p}}\mathcal{L}(\boldsymbol{w}_{t+1},\boldsymbol{p}_{t+1}) -\boldsymbol{z}_{\boldsymbol{p}}^t-(\boldsymbol{z}_{\boldsymbol{p}}^{t+1}-\boldsymbol{z}_{\boldsymbol{p}}^t) \|_2^2]\nonumber\\
	=&\mathbb{E}[\| \nabla_{\boldsymbol{p}}\mathcal{L}(\boldsymbol{w}_{t+1},\boldsymbol{p}_{t+1}) -\boldsymbol{z}_{\boldsymbol{p}}^t+b\boldsymbol{z}_{\boldsymbol{p}}^t-bH^{t+1} \|_2^2]\nonumber\\
	=&\mathbb{E}[\|(1-b)\left(\nabla_{\boldsymbol{p}}\mathcal{L}(\boldsymbol{w}_t,\boldsymbol{p}_t)-\boldsymbol{z}_{\boldsymbol{p}}^t\right)+(1-b)\left( \nabla_{\boldsymbol{p}}\mathcal{L}(\boldsymbol{w}_{t+1},\boldsymbol{p}_{t+1})-\nabla_{\boldsymbol{p}}\mathcal{L}(\boldsymbol{w}_{t},\boldsymbol{p}_{t}) \right) +b\left( \nabla_{\boldsymbol{p}}\mathcal{L}(\boldsymbol{w}_{t+1},\boldsymbol{p}_{t+1})-H^{t+1} \right)\|_2^2]\nonumber\\
	=&(1-b)^2\mathbb{E}[\|\nabla_{\boldsymbol{p}}\mathcal{L}(\boldsymbol{w}_t,\boldsymbol{p}_t)-\boldsymbol{z}_{\boldsymbol{p}}^t\|_2^2]+(1-b)^2\mathbb{E}[\| \nabla_{\boldsymbol{p}}\mathcal{L}(\boldsymbol{w}_{t+1},\boldsymbol{p}_{t+1})-\nabla_{\boldsymbol{p}}\mathcal{L}(\boldsymbol{w}_{t},\boldsymbol{p}_{t}) \|_2^2] \nonumber\\
	&+b^2\mathbb{E}[\| \nabla_{\boldsymbol{p}}\mathcal{L}(\boldsymbol{w}_{t+1},\boldsymbol{p}_{t+1})-H^{t+1} \|_2^2]+(1-b)^2\langle \nabla_{\boldsymbol{p}}\mathcal{L}(\boldsymbol{w}_t,\boldsymbol{p}_t)-\boldsymbol{z}_{\boldsymbol{p}}^t, \nabla_{\boldsymbol{p}}\mathcal{L}(\boldsymbol{w}_{t+1},\boldsymbol{p}_{t+1})-\nabla_{\boldsymbol{p}}\mathcal{L}(\boldsymbol{w}_{t},\boldsymbol{p}_{t})\rangle\nonumber\\
	=&(1-b)^2\mathbb{E}[\|\nabla_{\boldsymbol{p}}\mathcal{L}(\boldsymbol{w}_t,\boldsymbol{p}_t)-\boldsymbol{z}_{\boldsymbol{p}}^t\|_2^2]+(1-b)^2\mathbb{E}[\| \nabla_{\boldsymbol{p}}\mathcal{L}(\boldsymbol{w}_{t+1},\boldsymbol{p}_{t+1})-\nabla_{\boldsymbol{p}}\mathcal{L}(\boldsymbol{w}_{t},\boldsymbol{p}_{t}) \|_2^2] \nonumber\\
	&+b^2\mathbb{E}[\| \nabla_{\boldsymbol{p}}\mathcal{L}(\boldsymbol{w}_{t+1},\boldsymbol{p}_{t+1})-H^{t+1} \|_2^2]+(1-b)^2\mathbb{E}[\langle \nabla_{\boldsymbol{p}}\mathcal{L}(\boldsymbol{w}_t,\boldsymbol{p}_t)-\boldsymbol{z}_{\boldsymbol{p}}^t, \nabla_{\boldsymbol{p}}\mathcal{L}(\boldsymbol{w}_{t+1},\boldsymbol{p}_{t+1})-\nabla_{\boldsymbol{p}}\mathcal{L}(\boldsymbol{w}_{t},\boldsymbol{p}_{t})\rangle]\nonumber\\
	\leq&(1-b)^2(1+b)\mathbb{E}[\|\nabla_{\boldsymbol{p}}\mathcal{L}(\boldsymbol{w}_t,\boldsymbol{p}_t)-\boldsymbol{z}_{\boldsymbol{p}}^t\|_2^2]+(1-b)^2(1+\dfrac{1}{b})\| \nabla_{\boldsymbol{p}}\mathcal{L}(\boldsymbol{w}_{t+1},\boldsymbol{p}_{t+1})-\nabla_{\boldsymbol{p}}\mathcal{L}(\boldsymbol{w}_{t},\boldsymbol{p}_{t}) \|_2^2 \nonumber\\
	&+b^2\mathbb{E}[\| \nabla_{\boldsymbol{p}}\mathcal{L}(\boldsymbol{w}_{t+1},\boldsymbol{p}_{t+1})-H^{t+1} \|_2^2]\nonumber\\
	\leq&(1-b)\mathbb{E}[\|\nabla_{\boldsymbol{p}}\mathcal{L}(\boldsymbol{w}_t,\boldsymbol{p}_t)-\boldsymbol{z}_{\boldsymbol{p}}^t\|_2^2]+\dfrac{1}{b}L^2\mathbb{E}[\|\boldsymbol{w}_{t+1}-\boldsymbol{w}_t\|_2^2+\|\boldsymbol{p}_{t+1}-\boldsymbol{p}_t\|_2^2]+b^2\sigma_2^2.
	\end{align}
	Similarly, we have 
	\begin{align}
	&\mathbb{E}[\| \nabla_{\boldsymbol{w}}\mathcal{L}(\boldsymbol{w}_{t+1},\boldsymbol{p}_{t+1}) -\boldsymbol{z}_{\boldsymbol{w}}^{t+1} \|_2^2]\nonumber\\
	\leq&(1-b)\mathbb{E}[\|\nabla_{\boldsymbol{w}}\mathcal{L}(\boldsymbol{w}_t,\boldsymbol{p}_t)-\boldsymbol{z}_{\boldsymbol{w}}^t\|_2^2]+\dfrac{1}{b}L^2\mathbb{E}[\|\boldsymbol{w}_{t+1}-\boldsymbol{w}_t\|_2^2+\|\boldsymbol{p}_{t+1}-\boldsymbol{p}_t\|_2^2]+b^2\sigma_1^2.
	\end{align}
	
\end{proof}

\subsection{Proof of Theorem \ref{theorem:convergence_rate_of_ASZOGD}}
\begin{proof}
	Summing up the inequality in Lemma \ref{lemma:bound_of_z_in_ASZOGD}, we have
	\begin{align}
	&\sum_{t=1}^{T}\mathbb{E}[\| \nabla_{\boldsymbol{w}}\mathcal{L}(\boldsymbol{w}_{t+1},\boldsymbol{p}_{t+1}) -\boldsymbol{z}_{\boldsymbol{w}}^{t+1} \|_2^2]\nonumber\\
	\leq&(1-b)\sum_{t=1}^{T}\mathbb{E}[\|\nabla_{\boldsymbol{w}}\mathcal{L}(\boldsymbol{w}_t,\boldsymbol{p}_t)-\boldsymbol{z}_{\boldsymbol{w}}^t\|_2^2]+\dfrac{1}{b}L^2\sum_{t=1}^{T}\mathbb{E}[\|\boldsymbol{w}_{t+1}-\boldsymbol{w}_t\|_2^2+\|\boldsymbol{p}_{t+1}-\boldsymbol{p}_t\|_2^2]+b^2\sigma_1^2T.
	\end{align}
	Then, we have 
	\begin{align}
	&\sum_{t=1}^{T}\mathbb{E}[\| \nabla_{\boldsymbol{w}}\mathcal{L}(\boldsymbol{w}_{t+1},\boldsymbol{p}_{t+1}) -\boldsymbol{z}_{\boldsymbol{w}}^{t+1} \|_2^2]\nonumber\\
	\leq&\dfrac{1}{b}\mathbb{E}[\|\nabla_{\boldsymbol{w}}\mathcal{L}(\boldsymbol{w}_1,\boldsymbol{p}_1)-\boldsymbol{z}_{\boldsymbol{w}}^1\|_2^2]+\dfrac{1}{b^2}L^2\sum_{t=1}^{T}\mathbb{E}[\|\boldsymbol{w}_{t+1}-\boldsymbol{w}_t\|_2^2+\|\boldsymbol{p}_{t+1}-\boldsymbol{p}_t\|_2^2]+b\sigma_1^2T.
	\end{align}
	Similarly, we have 
	\begin{align}
	&\sum_{t=1}^{T}\mathbb{E}[\| \nabla_{\boldsymbol{p}}\mathcal{L}(\boldsymbol{w}_{t+1},\boldsymbol{p}_{t+1}) -\boldsymbol{z}_{\boldsymbol{p}}^{t+1} \|_2^2]\nonumber\\
	\leq&\dfrac{1}{b}\mathbb{E}[\|\nabla_{\boldsymbol{w}}\mathcal{L}(\boldsymbol{w}_1,\boldsymbol{p}_1)-\boldsymbol{z}_{\boldsymbol{w}}^1\|_2^2]+\dfrac{1}{b^2}L^2\sum_{t=1}^{T}\mathbb{E}[\|\boldsymbol{w}_{t+1}-\boldsymbol{w}_t\|_2^2+\|\boldsymbol{p}_{t+1}-\boldsymbol{p}_t\|_2^2]+b\sigma_1^2T.
	\end{align}
	
	\begin{align}
	&\sum_{t=1}^{T }\|\boldsymbol{p}_{t+1}-\boldsymbol{p}^*(\boldsymbol{w}_{t+1})\|_2^2\nonumber\\
	\leq&\dfrac{4}{\tau a \eta_{p}c_{2,l}}(\|\boldsymbol{p}^*(\boldsymbol{w}_{1})-\boldsymbol{p}_1\|_2^2-\dfrac{1}{4a}\sum_{t=1}^{T}\|\boldsymbol{p}_t-{\boldsymbol{p}}_{t+1}\|_2^2+\dfrac{8a\eta_pc_{2,l}}{\tau}\sum_{t=1}^{T}\|\nabla_{\boldsymbol{p}}\mathcal{L}(\boldsymbol{w}_t,\boldsymbol{p}_t)-\boldsymbol{z}_{\boldsymbol{p}}^t\|_2^2\nonumber\\
	&+\dfrac{8L^2_g}{\tau a \eta_{p}c_{2,l}}\sum_{t=1}^{T}\|\boldsymbol{w}_{t}-\boldsymbol{w}_{t+1}\|_2^2)\nonumber\\
	\leq&\dfrac{4}{\tau a \eta_{p}c_{2,l}}\|\boldsymbol{p}^*(\boldsymbol{w}_{1})-\boldsymbol{p}_1\|_2^2-\dfrac{1}{\tau a^2 \eta_{p}c_{2,l}}\sum_{t=1}^{T}\|\boldsymbol{p}_t-{\boldsymbol{p}}_{t+1}\|_2^2+\dfrac{32}{\tau^2}\sum_{t=1}^{T}\|\nabla_{\boldsymbol{p}}\mathcal{L}(\boldsymbol{w}_t,\boldsymbol{p}_t)-\boldsymbol{z}_{\boldsymbol{p}}^t\|_2^2+\dfrac{32L^2_g\eta_{w}^2c_{1,u}^2}{\tau^2 a^2 \eta_{p}^2c_{2,l}^2}\sum_{t=1}^{T}\|\boldsymbol{z}_{\boldsymbol{w}}^{t        }\|_2^2. 
	\end{align}
	Thus, we have 
	\begin{align}
	&\mathbb{E}[\sum_{t=1}^{T }\|\boldsymbol{p}_{t+1}-\boldsymbol{p}^*(\boldsymbol{w}_{t+1})\|_2^2]\nonumber\\
	\leq&\dfrac{4}{\tau a \eta_{p}c_{2,l}}\mathbb{E}[\|\boldsymbol{p}^*(\boldsymbol{w}_{1})-\boldsymbol{p}_1\|_2^2]+\mathbb{E}[\dfrac{32L^2_g\eta_{w}^2c_{1,u}^2}{\tau^2 a^2 \eta_{p}^2c_{2,l}^2}\sum_{t=1}^{T}\|\boldsymbol{z}_{\boldsymbol{w}}^{t        }\|_2^2 -\dfrac{1}{\tau a^2 \eta_{p}c_{2,l}}\sum_{t=1}^{T}\|\boldsymbol{p}_t-{\boldsymbol{p}}_{t+1}\|_2^2] \nonumber\\
	&+\dfrac{32}{\tau^2b}\mathbb{E}[\|\nabla_{\boldsymbol{p}}\mathcal{L}(\boldsymbol{w}_1,\boldsymbol{p}_1)-\boldsymbol{z}_{\boldsymbol{p}}^1\|_2^2]+\mathbb{E}\left[ \sum_{t=1}^{T}\dfrac{32}{\tau^2}\left( b\sigma_2^2+\dfrac{L^2(\|\boldsymbol{w}_{t+1}-\boldsymbol{w}_t\|_2^2+\|\boldsymbol{p}_{t+1}-\boldsymbol{p}_t\|_2^2)}{b^2} \right) \right]\nonumber\\
	\leq&\dfrac{4}{\tau a \eta_{p}c_{2,l}}\mathbb{E}[\|\boldsymbol{p}^*(\boldsymbol{w}_{1})-\boldsymbol{p}_1\|_2^2]+\mathbb{E}\left[\left(\dfrac{32L^2_g\eta_{w}^2c_{1,u}^2}{\tau^2 a^2 \eta_{p}^2c_{2,l}^2}+\dfrac{32L^2\eta_{w}^2c_{1,u}^2}{\tau^2b^2}\right)\sum_{t=1}^{T}\|\boldsymbol{z}_{\boldsymbol{w}}^{t        }\|_2^2\right] \nonumber\\
	&+\dfrac{32}{\tau^2b}\mathbb{E}[\|\nabla_{\boldsymbol{p}}\mathcal{L}(\boldsymbol{w}_1,\boldsymbol{p}_1)-\boldsymbol{z}_{\boldsymbol{p}}^1\|_2^2]+\mathbb{E}\left[ \sum_{t=1}^{T}\dfrac{32}{\tau^2} b\sigma_2^2 \right]\nonumber\\
	&+\mathbb{E}\left[\left(\dfrac{32L^2}{\tau^2b^2}-\dfrac{1}{\tau a^2 \eta_{p}c_{2,l}}\right)\sum_{t=1}^{T}\|\boldsymbol{p}_t-{\boldsymbol{p}}_{t+1}\|_2^2] \right].
	\end{align}
	In addition, we have 
	\begin{align}
	&\dfrac{1}{2}\|\nabla g(\boldsymbol{w}_t)\|_2^2\nonumber\\
	\leq&\dfrac{g(\boldsymbol{w}_t)-g(\boldsymbol{w}_{t+1})}{\eta_wc_{1,l}}-\dfrac{1}{4}\|\boldsymbol{z}_{\boldsymbol{w}}^t\|_2^2+\dfrac{\mu^2 L^2(d+3)^3}{2}+L^2\|\boldsymbol{p}^*(\boldsymbol{w}_t)-\boldsymbol{p}_t\|_2^2+\|\nabla_{\boldsymbol{w}}\mathcal{L}(\boldsymbol{w}_t,\boldsymbol{p}_t)-\boldsymbol{z}_{\boldsymbol{w}}^t\|_2^2.
	\end{align}
	
	Summing up from $t=1,\cdots,T$ and taking expectation, we have 
	\begin{align}
	&\mathbb{E}[\dfrac{1}{2}\sum_{t=1}^{T}\|\nabla g(\boldsymbol{w}_t)\|_2^2]\nonumber\\
	\leq&\dfrac{g(\boldsymbol{w}_1)-g(\boldsymbol{w}_{T+1})}{\eta_wc_{1,l}}-\dfrac{1}{4}\mathbb{E}[\sum_{t=1}^{T}\|\boldsymbol{z}_{\boldsymbol{w}}^t\|_2^2]+\dfrac{\mu^2 TL^2(d+3)^3}{2}\nonumber\\
	&+\dfrac{4L^2}{\tau a \eta_{p}c_{2,l}}\mathbb{E}[\|\boldsymbol{p}^*(\boldsymbol{w}_{1})-\boldsymbol{p}_1\|_2^2]+\mathbb{E}\left[\left(\dfrac{L^2\eta_{w}^2c_{1,u}^2}{b^2}+\dfrac{32L^2_gL^2\eta_{w}^2c_{1,u}^2}{\tau^2 a^2 \eta_{p}^2c_{2,l}^2}+\dfrac{32L^4\eta_{w}^2c_{1,u}^2}{\tau^2b^2}\right)\sum_{t=1}^{T}\|\boldsymbol{z}_{\boldsymbol{w}}^{t        }\|_2^2\right] \nonumber\\
	&+\dfrac{32L^2}{\tau^2b}\mathbb{E}[\|\nabla_{\boldsymbol{p}}\mathcal{L}(\boldsymbol{w}_1,\boldsymbol{p}_1)-\boldsymbol{z}_{\boldsymbol{p}}^1\|_2^2]+\mathbb{E}\left[ \sum_{t=1}^{T}\dfrac{32L^2}{\tau^2} b\sigma_2^2 \right]\nonumber\\
	&+\mathbb{E}\left[\left(\dfrac{L^2}{b^2}+\dfrac{32L^4}{\tau^2b^2}-\dfrac{L^2}{\tau a^2 \eta_{p}c_{2,l}}\right)\sum_{t=1}^{T}\|\boldsymbol{p}_t-{\boldsymbol{p}}_{t+1}\|_2^2] \right]\nonumber\\
	&+\dfrac{1}{b}\mathbb{E}[\|\nabla_{\boldsymbol{w}}\mathcal{L}(\boldsymbol{w}_1,\boldsymbol{p}_1)-\boldsymbol{z}_{\boldsymbol{w}}^1\|_2^2]+b\sigma_1^2T.
	\end{align}
	Let $\boldsymbol{p}^*(w_{1})=\boldsymbol{p}_1$, $\boldsymbol{z}_{\boldsymbol{p}}^1=H(\boldsymbol{w}_{t},\boldsymbol{p}_{t},f_{\mathcal{M}_3})$, $\boldsymbol{z}_{\boldsymbol{w}}^1=G_{\mu}^{\mathcal{L}}(\boldsymbol{w}_{t},\boldsymbol{p}_{t},\ell_{\mathcal{M}_1},f_{\mathcal{M}_2},\boldsymbol{u}_{[q]})$, $\eta_{p}\leq\min\{\dfrac{b^2}{\tau a^2c_{2,l}},\dfrac{\tau b^2}{32L^2a^2c_{2,l}}\}$, $\eta_{w}^2\leq\min\{\dfrac{b^2}{4c_{1,u}^2L^2},\dfrac{\tau^2a^2\eta_{p}^2c_{2,l}^2}{128L_g^2L^2c_{1,u}},\dfrac{\tau^2b^2}{128L^4c_{1,u}^2}\}$, we have 
	\begin{align}
	&\dfrac{1}{T}\mathbb{E}[\sum_{t=1}^{T}\|\nabla g(\boldsymbol{w}_t)\|_2^2]
	\leq\dfrac{2(g(\boldsymbol{w}_1)-g(\boldsymbol{w}_{T+1}))}{T\eta_wc_{1,l}}+\dfrac{64L^2}{\tau^2 bT}\sigma_2^2+\dfrac{2\sigma_1^2}{bT}+\mu^2 L^2(d+3)^3+\dfrac{64L^2}{\tau^2} b\sigma_2^2 +2b\sigma_1^2.
	\end{align}
	Bound the left term by $\epsilon^2$, we have $\mu\leq\dfrac{\epsilon}{L(d+3)^{3/2}}$, $b\leq\min\{\dfrac{\epsilon^2}{2\sigma_1^2},\dfrac{\tau^2\epsilon^2}{64\sigma_2^2L^2}\}$ and $T\geq\max\{\dfrac{2(g(\boldsymbol{w}_1)-g(\boldsymbol{w}_T))}{\epsilon^2\eta_{w}c_{1,l}},\dfrac{2\sigma_1^2}{\epsilon^2b},\dfrac{64\sigma_2^2L^2}{\epsilon^2\tau^2 b}\}$.
	
\end{proof}
\section{Additional Experiments}
\subsection{Impact of the Hyper-parameters}
In this section, we discuss the impact of the learning rate $\eta_{\boldsymbol{p}}$ and $\eta_{\boldsymbol{w}}$. First, let $\eta_{\boldsymbol{p}}=0.001$. We evaluate the performance of our method in two applications with $\eta_{\boldsymbol{w}}$ chosen from $\{0.01,0.001,0.0001\}$. Then, let $\eta_{\boldsymbol{w}}=0.01$. We evaluate the performance in two applications with $\eta_{\boldsymbol{p}}$ chosen from $\{0.01,0.001,0.0001\}$. All the experiments are presented in Tables \ref{tab:fairness1}, \ref{tab:fairness2}, \ref{tab:pairwise1} and \ref{tab:pairwise2}. We can find that it is important to choose a proper $\eta_{\boldsymbol{w}}$. In addition, our method is not sensitive to $\eta_{\boldsymbol{p}}$. We also give the results of our method with different $\beta$ $a$ and $b$ in Table \ref{tab:nonlinear_beta} and Table \ref{tab:nonlinear_ab}. We can find that our method is not sensitive to $\beta$, $a$ and $b$.

\subsection{Performance with Nonlinear Model}
\begin{figure}[t]
	\centering
	\small
	\begin{subfigure}[b]{0.245\textwidth}
		\centering
		\includegraphics[width=1.8in]{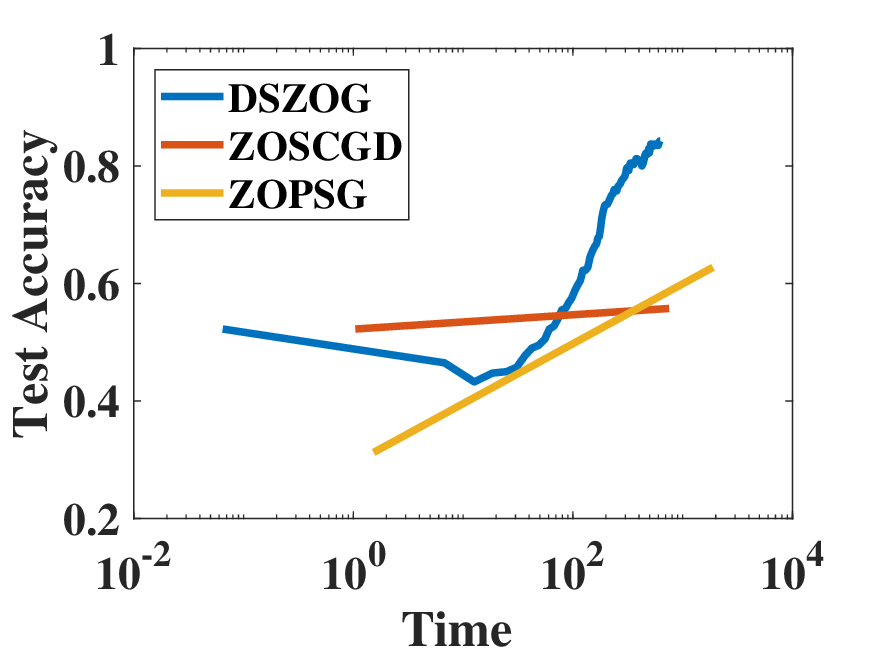}
		\caption{D1}
	\end{subfigure}
	\begin{subfigure}[b]{0.245\textwidth}
		\centering
		\includegraphics[width=1.8in]{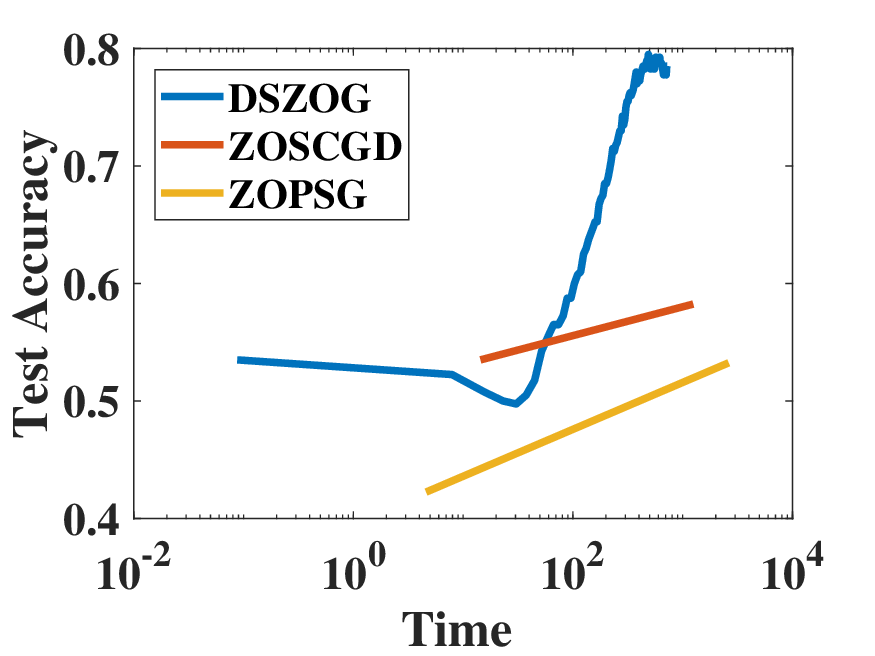}
		\caption{D2}
	\end{subfigure}
	\begin{subfigure}[b]{0.245\textwidth}
		\centering
		\includegraphics[width=1.8in]{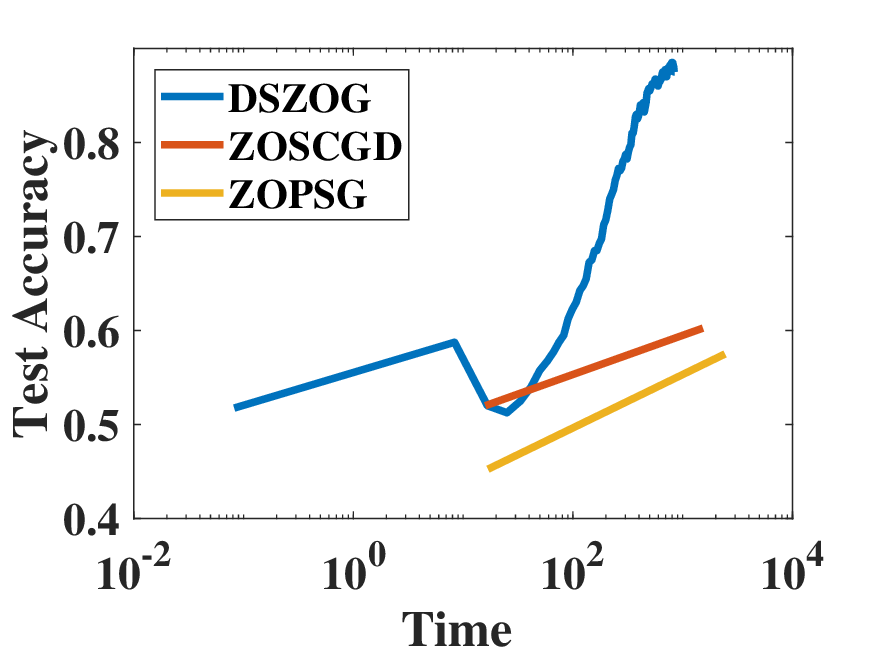}
		\caption{D3}
	\end{subfigure}	
	\begin{subfigure}[b]{0.245\textwidth}
		\centering
		\includegraphics[width=1.8in]{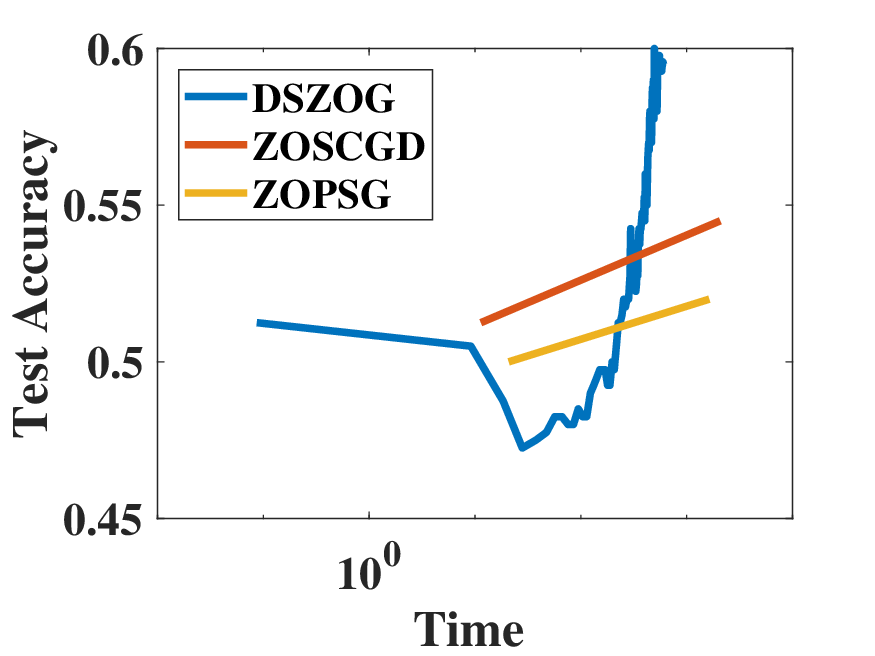}
		\caption{D4}
		
	\end{subfigure}
	\caption{Performance of our method in fairness with kernel method.}
	\label{fig:fairness2}
\end{figure}

We use the kernel method $k(x,x')=\exp(-\gamma\|x-x'\|_2^2)$ to conduct a nonlinear model. The hyperparameter of all the methods are set according to Section 6. The results of fairness are presented in Table \ref{tab:nonlinear} and Figure \ref{fig:fairness2}. We can find that our method is still superior than other methods.

\begin{table*}[]
	\centering
	\caption{Test accuracy (\%) of DSZOG in classification with fairness constraints when using different $\eta_{\boldsymbol{w}}$.}
	\label{tab:fairness1}
	\begin{tabular}{llll}
		\toprule
		Data& 0.01                & 0.001        & 0.0001                   \\
		\hline
		D1    & $87.33\pm0.38$   & $87.21\pm0.34$ & $80.68\pm0.32$  \\
		\hline
		D2 & $84.75\pm0.25$  & $83.01\pm0.24$ & $67.34\pm0.78$  \\
		\hline
		D3 & $83.58\pm0.14$ & $83.78\pm0.45$ & $77.22\pm0.32$\\
		\hline
		D4   & $64.91\pm0.94$  & $63.27\pm0.45$ &$55.13\pm0.46$ \\
		\bottomrule
	\end{tabular}
	
\end{table*}
\begin{table*}[]
	\centering
	\caption{Test accuracy (\%) of DSZOG in classification with pairwise constraints when using different $\eta_{\boldsymbol{w}}$.}
	\label{tab:pairwise1}
	\begin{tabular}{llll}
		\toprule
		Data& 0.01                & 0.001        & 0.0001                   \\
		\hline
		a9a    & $75.90\pm0.26$   & $74.45\pm0.67$ & $65.33\pm0.43$  \\
		\hline
		w8a & $89.94\pm0.28$  & $88.28\pm0.80$ & $56.88\pm0.22$  \\
		\hline
		gen & $82.33\pm0.76$  & $81.78\pm0.45$ & $57.34\pm0.22$\\
		\hline
		svm   & $79.56\pm0.49$  & $79.56\pm0.45$ & $56.23\pm0.75$\\
		\bottomrule
	\end{tabular}
	
\end{table*}

\begin{table*}[]
	\centering
	\caption{Test accuracy (\%) of DSZOG in classification with fairness constraints when using different $\eta_{\boldsymbol{p}}$.}
	\label{tab:fairness2}
	\begin{tabular}{llll}
		\toprule
		Data& 0.01                & 0.001        & 0.0001                   \\
		\hline
		D1    & $87.33\pm0.38$   & $86.12\pm0.23$ & $85.13\pm0.22$  \\
		\hline
		D2 & $84.75\pm0.25$  & $83.22\pm0.14$ & $83.05\pm0.03$  \\
		\hline
		D3 & $83.58\pm0.14$ & $83.29\pm0.31$ & $81.18\pm0.22$\\
		\hline
		D4   & $64.91\pm0.94$  & $62.12\pm0.31$ &$62.23\pm0.22$ \\
		\bottomrule
	\end{tabular}
	
\end{table*}
\begin{table*}[]
	\centering
	\caption{Test accuracy (\%) of DSZOG in classification with pairwise constraints when using different $\eta_{\boldsymbol{p}}$.}
	\label{tab:pairwise2}
	\begin{tabular}{llll}
		\toprule
		Data& 0.01                & 0.001        & 0.0001                   \\
		\hline
		a9a    & $75.90\pm0.26$   & $73.33\pm0.31$ & $73.64\pm0.22$  \\
		\hline
		w8a & $89.94\pm0.28$  & $89.23\pm0.17$ & $88.82\pm0.56$  \\
		\hline
		gen & $82.33\pm0.76$  & $82.11\pm0.64$ & $81.28\pm0.34$\\
		\hline
		svm   & $79.56\pm0.49$  & $78.56\pm0.44$ & $78.45\pm0.76$\\
		\bottomrule
	\end{tabular}
	
\end{table*}

\begin{table}[]
	\centering
	\small
	\setlength{\tabcolsep}{1mm}
	\caption{Test accuracy (\%) in fairness with kernel method.}\label{tab:nonlinear}
	\begin{tabular}{llll}
		\toprule
		Data& DSZOG                & ZOSCGD        & ZOPSGD                   \\
		\hline
		D1    & $\textbf{84.25}$   & $55.75$ & $62.75$  \\
		\hline
		D2 & $\textbf{78.50}$  & $58.25$ & $53.25$ \\
		\hline
		D3 & $\textbf{87.50}$ & $60.25$ & $57.50$  \\
		\hline
		D4   & $\textbf{59.52}$  & $54.50$ & $52.25$ \\
		\bottomrule
	\end{tabular}
\end{table}
\begin{table}[]
	\centering
	\small
	\setlength{\tabcolsep}{1mm}
	\caption{Test accuracy of DSZOG with different $\beta$ fairness ($\eta_w=0.001$, $\eta_p=0.1$, $a=0.9$, $b=0.9$).}\label{tab:nonlinear_beta}
	\begin{tabular}{llll}
		\toprule
		Data& 0.1                & 1        & 10                   \\
		\hline
		D1    & ${84.25}$   & $84.73$ & $83.12$  \\
		\hline
		D2 & ${78.50}$  & $77.98$ & $78.33$ \\
		\hline
		D3 & ${87.50}$ & $87.67$ & $87.85$  \\
		\hline
		D4   & ${59.52}$  & $58.61$ & $58.98$ \\
		\bottomrule
	\end{tabular}
\end{table}
\begin{table}[]
	\centering
	\small
	\setlength{\tabcolsep}{1mm}
	\caption{Test accuracy (\%) of DSZOG with different $a$ and $b$ fairness constraints ($\eta_w=0.001$, $\eta_p=0.1$, $\beta=0.1$).}\label{tab:nonlinear_ab}
\begin{tabular}{l|lll|lll}
\hline
   & \multicolumn{3}{l|}{$a=0.9$}                                          & \multicolumn{3}{l}{$b=0.9$}                                           \\ \hline
   & \multicolumn{1}{l|}{$b=0.1$} & \multicolumn{1}{l|}{$b=0.5$} & $b=0.9$ & \multicolumn{1}{l|}{$a=0.1$} & \multicolumn{1}{l|}{$a=0.5$} & $b=0.9$ \\ \hline
D1 & \multicolumn{1}{l|}{$83.93$}        & \multicolumn{1}{l|}{$84.12$}        & $84.25$        & \multicolumn{1}{l|}{$82.12$}        & \multicolumn{1}{l|}{$83.21$}        &   $84.25$      \\ \hline
D2 & \multicolumn{1}{l|}{$78.33$}        & \multicolumn{1}{l|}{$78.27$}        & $78.50$        & \multicolumn{1}{l|}{$77.29$}        & \multicolumn{1}{l|}{$78.45$}        &  $78.50$       \\ \hline
D3 & \multicolumn{1}{l|}{$85.34$}        & \multicolumn{1}{l|}{$87.410$}        &$87.50$         & \multicolumn{1}{l|}{$87.23$}        & \multicolumn{1}{l|}{$87.14$}        & $87.50$        \\ \hline
D4 & \multicolumn{1}{l|}{$58.82$}        & \multicolumn{1}{l|}{$57.23$}        & $59.52$        & \multicolumn{1}{l|}{$58.78$}        & \multicolumn{1}{l|}{$59.48$}        &    $59.52$     \\ \hline
\end{tabular}
\end{table}

\end{document}